\documentclass[12pt]{amsart}
\usepackage{amsbsy}
\usepackage{graphicx,epsfig,subfigure,psfrag}
\begin{document}

\newtheorem{theorem}{Theorem}
\newtheorem{proposition}{Proposition}
\newtheorem{lemma}{Lemma}
\newtheorem{corollary}{Corollary}
\newtheorem{definition}{Definition}
\newtheorem{remark}{Remark}
\newcommand{\be}{\begin{equation}}
\newcommand{\ee}{\end{equation}}
\newcommand{\tex}{\textstyle}
\numberwithin{equation}{section} \numberwithin{theorem}{section}
\numberwithin{proposition}{section} \numberwithin{lemma}{section}
\numberwithin{corollary}{section}
\numberwithin{definition}{section} \numberwithin{remark}{section}
\newcommand{\ren}{\mathbb{R}^N}
\newcommand{\re}{\mathbb{R}}
\newcommand{\n}{\nabla}
\newcommand{\iy}{\infty}
\newcommand{\pa}{\partial}
\newcommand{\fp}{\noindent}
\newcommand{\ms}{\medskip\vskip-.1cm}
\newcommand{\mpb}{\medskip}
\newcommand{\BB}{{\bf B}}
\newcommand{\AAA}{{\bf A}}
\newcommand{\Am}{{\bf A}_{2m}}
\newcommand{\ef}{\eqref}
\newcommand{\eee}{{\mathrm e}}
\newcommand{\ii}{{\mathrm i}}
\renewcommand{\a}{\alpha}
\renewcommand{\b}{\beta}
\newcommand{\g}{\gamma}
\newcommand{\G}{\Gamma}
\renewcommand{\d}{\delta}
\newcommand{\D}{\Delta}
\newcommand{\e}{\varepsilon}
\newcommand{\var}{\varphi}
\renewcommand{\l}{\lambda}
\renewcommand{\o}{\omega}
\renewcommand{\O}{\Omega}
\newcommand{\s}{\sigma}
\renewcommand{\t}{\tau}
\renewcommand{\th}{\theta}
\newcommand{\z}{\zeta}
\newcommand{\wx}{\widetilde x}
\newcommand{\wt}{\widetilde t}
\newcommand{\noi}{\noindent}
\newcommand{\uu}{{\bf u}}
\newcommand{\UU}{{\bf U}}
\newcommand{\VV}{{\bf V}}
\newcommand{\ww}{{\bf w}}
\newcommand{\vv}{{\bf v}}
\newcommand{\WW}{{\bf W}}
\newcommand{\hh}{{\bf h}}
\newcommand{\di}{{\rm div}\,}
\newcommand{\inA}{\quad \mbox{in} \quad \ren \times \re_+}
\newcommand{\inB}{\quad \mbox{in} \quad}
\newcommand{\inC}{\quad \mbox{in} \quad \re \times \re_+}
\newcommand{\inD}{\quad \mbox{in} \quad \re}
\newcommand{\forA}{\quad \mbox{for} \quad}
\newcommand{\whereA}{,\quad \mbox{where} \quad}
\newcommand{\asA}{\quad \mbox{as} \quad}
\newcommand{\andA}{\quad \mbox{and} \quad}
\newcommand{\withA}{,\quad \mbox{with} \quad}
\newcommand{\orA}{,\quad \mbox{or} \quad}
\newcommand{\ssk}{\smallskip}
\newcommand{\LongA}{\quad \Longrightarrow \quad}
\def\com#1{\fbox{\parbox{6in}{\texttt{#1}}}}
\def\N{{\mathbb N}}
\def\A{{\cal A}}
\def\WW{{\cal W}}
\newcommand{\de}{\,d}
\newcommand{\eps}{\varepsilon}
\newcommand{\spt}{{\mbox spt}}
\newcommand{\ind}{{\mbox ind}}
\newcommand{\supp}{{\mbox supp}}
\newcommand{\dip}{\displaystyle}
\newcommand{\prt}{\partial}
\renewcommand{\theequation}{\thesection.\arabic{equation}}
\renewcommand{\baselinestretch}{1.2}

\title
{\bf Towards the KPP--problem  and
  ${\bf { log\, t}}$--front \\ shift for higher-order nonlinear PDEs
II.\\ Quasilinear bi- and tri-harmonic  equations}

\author{
V.A.~Galaktionov}

\address{Department of Mathematical Sciences, University of Bath,
 Bath BA2 7AY, UK}
\email{masvg@bath.ac.uk}



  \keywords{KPP--problem, travelling waves, stability,  higher-order quasilinear parabolic
  equations,  $\log t$-front shifting}
 \subjclass{35K55, 35K40, 35K65}
 \date{\today}




\begin{abstract}

 Extensions of ideas of Kolmogorov, Petrovskii, and
Piskunov (1937) \cite{KPP} on travelling wave propagation in the
reaction-diffusion equation
 $$
 u_t=u_{xx}+u(1-u) \inB \re \times \re_+, \quad u_0(x)=H(-x) \equiv\{1\,\,
 \mbox{for}\,\, x<0; \,\,\, 0 \,\, \mbox{for} \,\, x \ge 0\},
 $$
 $H(\cdot)$ being the Heaviside function, are discussed. The
 present paper continues the study began in \cite{GKPPI}
 for higher-order
semilinear
  {\em bi-harmonic} and {\em tri-harmonic equations}
 $$
 u_t= -u_{xxxx} +u(1-u) \andA
u_t=u_{xxxxxx}+u(1-u),
 \quad \mbox{etc.}
 $$

Here,  some of the results are extended to the corresponding {\em
quasilinear} degenerate parabolic models with nonlinearities of
the {\em porous medium type} ($n>0$),
 $$
 u_t=-(|u|^n u)_{xxxx}+u(1-u), \quad u_t=(|u|^n u)_{xxxxxx}+u(1-u),
 \quad \mbox{etc.}
 $$

 Two main questions to discuss are:

 (i) existence of travelling waves via any analytical/numerical methods, and

 (ii) their stability and
 derivation of the
 $\log t$-shifting of moving fronts for some class of data $u_0$ (not for $H(-x)$).

\end{abstract}

\maketitle



\setcounter{equation}{0}
\section{Introduction: the classic KPP--problem, its quasilinear extensions, and higher-order
 quasilinear parabolic PDEs}
 \label{Sect1}
  \setcounter{equation}{0}











We begin with an introduction concerning classic results; see more
details in \cite[\S~1]{GKPPI}.

\subsection{The classic KPP--problem of 1937: convergence to TWs}

In the KPP--problem \cite{KPP} (1937)
 \be
 \label{1.1}
  u_t = u_{xx} +u(1-u), \quad x\in \re,\,\,
t>0; \quad u(x,0)=u_0(x) \,\,\, \mbox{in $\re$},
 \ee
    with
the step (Heaviside) initial function
 \be
 \label{1.H}
u_0(x) =H(-x) \equiv \left\{ \begin{matrix}1, \quad x<0;\\ 0,
\quad x\ge 0,
\end{matrix}
\right.
 \ee
the solution was proved to converge to the so-called {\em minimal}
travelling wave (TW) solution corresponding to the minimally
possible TW speed
 \be
 \label{minSp1}
 \mbox{the minimal speed of TW propagation}: \quad \l_0=2.
  \ee
  Looking for a TW profile $f(y)$ for
arbitrary $\l>0$ yields:
 \be
 \label{TW1}
 \left\{
  \begin{matrix}
 u_*(x,t)=f(y), \,\, y=x-\l t, \,\,\,\,\mbox{where $f$ solves the ODE problem}\ssk\ssk\ssk\\
  -\l f'=f''+f(1-f), \,\, y \in \re; \quad
 f(-\iy)=1, \,\, f(+\iy)=0. \qquad
 \end{matrix}
 \right.
  \ee
This 2nd-order ODE, on the phase-plane $\{f,f'\}$, by setting
$f'=P(f)$, reduces its order:
 \be
 \label{TW11}
  \tex{
  \frac {{\mathrm d}P}{{\mathrm d}f}=-\l- \frac {f(1-f)}P \, ,
}
 \ee
 and it was shown that there exists the {\em minimal} speed
 $\l_0=2$ and
the unique (up to translation)
 {\em minimal} TW profile $f(y)$.
Using the natural normalization
 \be
 \label {norm1}
 \tex{
 f(0)= \frac 12
 }
 \ee
 this minimal TW profile is defined uniquely.
  In addition,
  \be
  \label{norm2}
  f'(y) <0 \inB \re.
   \ee
The characteristic equation for the linearized operator in
\ef{TW1} has a multiple zero:
 \be
 \label{Norm21}
g'' + 2 g'+g=0 \andA g=\eee^{\mu y} \LongA (\mu+1)^2=0 \LongA
\mu_{1,2}=-1,
 \ee
 that yields the following asymptotic behaviour of $f(y)$:
 \be
  \label{norm23}
   f(y)=C_0 y \, {\mathrm e}^{-y}(1+o(1)) \asA y
  \to + \iy \whereA \mbox{$C_0>0$ is a constant.}
  \ee

 Concerning the relation between the
ODE TW problem \ef{TW11} for $\l_0=2$ and the PDE Cauchy one
\ef{1.1}, \ef{1.H}, the novel remarkable analysis in \cite{KPP} of
convergence as $t \to +\iy$ of the solution of the Cauchy problem
\ef{1.1}, \ef{1.H} to the minimal TW \ef{TW1} was performed in the
TW moving frame. This was very essential, and not in view of the
obvious $x$-translational invariance of the equation \ef{1.1}; see
below. Eventually, using PDE methods, the KPP-authors proved that
  the
TW front moves like
 \be
 \label{1.2}
  x_f(t) = 2t -
g(t) \quad \mbox{as $t\to \infty$}, \quad \mbox{with $g'(t) \to 0
$},
 \ee
 where the front location
 $x_f(t)$ is
uniquely determined from the equation
 \be
 \label{TW2}
  \tex{
 u(x_f(t),t) = \frac 12 \quad \mbox{for all $t \ge 0$}.
 }
  \ee
Then the convergence result of \cite{KPP} takes the form:
 \be
 \label{TW3}
 u(x_f(t)+y,t) \to f(y) \asA t \to + \iy \quad \mbox{uniformly in
 $y \in \re$}.
 \ee


\subsection{Bramson's $\log t$-front drift}

The next important question is the actual behaviour of the TW
shift $g(t)$ in \ef{1.2} for $t \gg 1$ (not addressed in
\cite{KPP}).
This is about an ``centre subspace (manifold) drift" of general
solutions of the KPP PDE \ef{1.1} along  a one-parameter family of
exact ODE TWs $\{f(y+a),\,\, a \in \re\}$.

\ssk

This open problem was solved in 1983
by Bramson \cite{Br} (see also \cite{Ga}) by using pure
probabilistic techniques (the Feynman--Kac integral formula
together with sample path estimates for Brownian motion).
It was proved that, within the PDE setting \ef{1.1}, \ef{1.H},
there is an {\em unbounded} $\log t$-shift of the moving TW front:
 \be
 \label{1.3}
  g(t) = k \log t(1+o(1)) \asA t \to +\iy, \quad
\mbox{with ${\bf k=  {\frac 32}}$},
  \ee
  Thus, \ef{1.3} implies eventual, as $t \to +\iy$, {\em infinite} retarding of the
  solution $u(x,t)$ from the corresponding minimal TW (uniquely fixed by \ef{norm1}), thought the
  convergence \ef{TW3} takes place in the TW frame.



\subsection{Known extensions: TW profiles for quasilinear
second-order reaction-diffusion equations}


The first attempts to extend the KPP results to quasilinear
second-order reaction-diffusion equations
 \be
 \label{RD33}
 u_t=(k(u)u_x)_x+Q(u), \quad Q(0)=Q(1)=0, \quad k(u) \ge 0,
  \ee
  including the porous medium diffusion operators with $k(u)=u^n$,
  $n>0$,
were known since the 1970s\footnote{It seems, the first ever
sufficiently general existence TW result for \ef{RD33} was
obtained by A.~Nepomnyashchy, which was published in Russian in
the technical Journal ``Voprosy Atomnoi Nauki i Tehniki"
(``Proceedings of Atomic Science and Technology") and was not
translated into English.}; see results and extra references in
\cite{And10, Atk81, Ar80, Fife77, GalTW81, Gild04, Gild05, Grin87,
Had75, Hos86, King03, Man2010, PaVa91}, \cite[p.~34]{SGKM},
\cite{San95, San97, Vaz07}, etc.

 However, it seems that any results concerning possible $\log t$-front shift  for quasilinear KPP problems
 including those with the degenerate
operators of the PME type
 \be
 \label{RD34}
 u_t=(u^{n+1})_{xx} + Q(u), \quad Q(0)=Q(1)=0, \quad n>0 \quad (u \ge
 0)
  \ee
  (with, say, positive non-Heaviside data; then $+\log t$-shift may occur)
  were still unknown.
Recall that, since \ef{RD34} describes finite propagation, for the
Heaviside initial data, the minimal travelling wave is indeed
stable (just by a standard barrier-comparison approach, this was
known since 1980s; see some references above). In this stability
sense, the quasilinear KPP-2 problem \ef{RD34} is simpler than the
classic one for $n=0$. However, we claim that $\log t$-shift is
possible for \ef{RD34} for some classes of initial data. Indeed,
by obvious reasons, in order to avoid  a trivial comparison from
below with a sufficiently shifted minimal TW $f(x-\l_0 t+a)$, with
$a \gg 1$, one needs that $u_0(x)$  must intersect {\sc all} such
TW profiles $f(x+a)$ (evidently, the Heaviside one does not do
that).
Actually, this means that, for \ef{RD34},  a
``correct" KPP-type setting (a really difficult one) assumes that
initial data {\em do not} have finite interface and are specially
distributed, that do not allow to use straightforward upper/lower
barriers for TW-like estimates above and below.


If fact, it is easy to construct a number of explicit TW solutions
for the PME with special source terms in \ef{RD34} using the {\em
pressure variable}:
 \be
 \label{pr1}
  \tex{
 u^n=v \LongA v_t=(n+1)vv_{xx}+ \frac{n+1}n\, (v_x)^2+q(v),
 \quad q(v)=nv^{\frac {n-1}n}\,Q(v^{\frac 1n}).
  }
  \ee
The ODE for TW profiles $f(y)$ then takes the form:
 \be
 \label{Pr2}
  \tex{
 -\l f'=(n+1)ff''+ \frac{n+1}n\, (f')^2 + q(f),
 }
  \ee
  and, as usual, we are looking for a solution $f$ having finite
  interface for $y>0$ and $f \to 1$ as $y \to -\iy$, with a linear
  spatial behaviour at the interfaces (this well corresponds to
  Darcy's law of their propagation for such weak solutions).

  \ssk

  \noi{\bf Example 1.} There exists the following explicit
  solution of \ef{Pr2} ($(\cdot)_+$ denotes the positive part):
   \be
   \label{Pr3}
    \begin{aligned}
   & f(y)=
    \tex{
    \frac {(-y)_+}{1-y}, \quad \l_0= \frac {n+1}n \whereA }
   \ssk\\
  &
 \tex{
   q(f)=\l(f-1)^2-2(n+1)(f-1)^3-(2(n+1)+ \frac{n+1}n)(f-1)^4.
   }
   \end{aligned}
   \ee

   \ssk

\noi{\bf Example 2.} There exists the following explicit
  solution of \ef{Pr2}:
   \be
   \label{Pr3N}
    \begin{aligned}
  & f(y)=
 \tex{
    \frac {(\eee^{-y}-1)_+}{\eee^{-y}+1}, \quad \l_0= \frac
    {n+1}{2n}}
   \whereA
    q(f)=
 \tex{
    \frac 12\,(1-f)[\l(1+f)
    }
 \ssk\\
     &
 \tex{
    -(n+1)f(1-f)  + (n+1)f(1+f)^2 - \frac 12\, \frac{n+1}n\,(1-f)(1+f)^2 ].
   }
    \end{aligned}
   \ee

   \ssk

We claim that, in the corresponding PDE framework,  such explicit
TWs admit a $\pm \log t$-front shifting along an ``affine centre
subspace", for suitable classes of initial data ($+\log t$ for
positive data, and $-\log t$ -- for changing sign one). The
derivations is similar (and simpler) to that in Section
\ref{SDiscr} for the quasilinear bi- and tri-harmonic flows.

\subsection{The present goal: extensions of the KPP--problem
to higher-order quasilinear parabolic PDEs}

The main goal of the present research comprising \cite{GKPPI,
GKPPIII}, and the present paper,  is to show that the
KPP--ideology
 can be extended to a variety of other more complicated
higher-order semilinear and quasilinear PDEs with source-type
terms.

Thus, our goal here  is to continue such a study, and, using
various analytic, formal, and  numerical methods), to show that
such a general viewing of the KPP ideas make sense, and that many
higher-order PDEs inherits some (but never all, of course) key
features of this classic analysis.

 Let us note that, in \cite{Bert01, Bert99} (see also a large list of  references therein, as well as those
  in the present paper),
 TW profiles were studied to a class of
 quasilinear {\em thin film equations}\footnote{Indeed, this quasilinear model is
 most relevant to the present study of \ef{PME4} and \ef{PME6}, though our solutions
 of the Cauchy problems
 {\em are not} nonnegative, are oscillatory, and
  of changing sign near the right-hand finite interface.}
  \be
  \label{TFE1}
  u_t+(f(u))_x= (b(u)u_x)_x-(c(u)u_{xxx})_x.
  \ee
  In particular, for more general TFE-type equations in $\ren \times \re_+$, stability of 1D TWs  was
  studied in \cite{Hu2008}, though, as in most of the papers
  mentioned around, the authors dealt with nonnegative solutions
  of a free-boundary value problem, rather than the Cauchy one.

\ssk

Thus, firstly, we will discuss some aspects of  KPP-type problems
for higher-order quasilinear partial differential equations
(PDEs), with the same Heaviside initial data. From applications,
such KPP-type problems deal with nonlinear higher-order diffusion
operators leading to well-known nowadays models; see references on
various fourth and $2m$th-order semilinear and quasilinear PDEs in
\cite{Bert01, Bert99, Collet90, Gl4, EGW1, Gal2m, GalCr}.

\ssk

To this end, in the present paper, we extend some of the KPP-type
results to the {\em quasilinear} bi-harmonic equation with the
diffusion operators of the fourth-order {\em porous medium} type
 \be
 \label{PME4}
  u_t= - (|u|^n u)_{xxxx}+ u(1-u) \inB \re \times \re_+ \whereA
  n>0.
   \ee
The corresponding TW with the speed of propagation $\l$ is then
governed by the following fourth-order ODE:
 \be
 \label{E5}
u_*(x,t)=f(y), \quad  y=x - \l t \LongA -\l f'=-(|f|^n f)
f^{(4)}+f(1-f),
 \ee
 with the singular boundary conditions at infinity:
 \be
 \label{BC1}
 f(y) =0 \,\,\,\mbox{for $y \gg 1$}  \andA f(y) \to 1 \asA y \to - \iy \quad
 \mbox{exponentially fast.}
  \ee
The first  condition in \ef{BC1} just takes into account that the
{\em degenerate} equation \ef{PME4} describes processes with {\em
finite propagation} of disturbances, a very well known fact for
such PDEs; see key references and results in \cite{GalRDE4n}.
Moreover, solutions are oscillatory and changing sign close to
finite interfaces. Therefore, the following example of a
nonnegative TW solution cannot be generic:

\ssk

\noi{\bf Example 3: nonnegative smooth TW for a bi-harmonic flow.}
We take $n=1$
 and consider the fourth-order parabolic equation with a TW $f \ge
 0$:
  \be
  \label{nn11}
  u_t=-(u^2)_{xxxx}+q(u) \LongA - \l f'=-(f^2)^{(4)}+q(f).
   \ee
   For a special choice of the reaction $q$ (see below), it admits the following explicit solution:
    \be
    \label{nn12}
     \tex{
  f(y)= \big[ \frac{(-y)}{1-y}\big]^3,
 }
   \ee
where the non-oscillatory behaviour $f(y) \sim (-y)^3$ close to
the interface $y=0^-$ is the actual decay for weak solutions  of
the Cauchy problem for $n=1$; see Section \ref{SectLocR}. It
follows from \ef{nn12} that
 \be
 \label{nn13}
  \tex{
  (-y)= \frac 1{1-f^{1/3}}-1,
  }
  \ee
  so that, substituting \ef{nn12} into the ODE in \ef{nn11},
  performing all the differentiations and eventually expressing
  $(-y)$ in terms of $f$ via \ef{nn13}, one obtains the actual
  source term $q(f)$, for which this is a solution. The required
  value  of the speed $\l$ is then obtained from
 the asymptotic analysis near the interface, as $y \to 0^-$, as
 follows:
  \be
  \label{nn15}
  f(y)=(-y)^3(1+O(y)) \LongA 3 \l y^2=-(y^6)^{(4)}+q(f)+...\, ,
  \,\,\, \mbox{so} \,\,\,
  \l_0=-120<0.
  \ee
  Indeed, for other values of $\l$, when the terms $O(y^2)$ are
  not cancelled in the equation in \ef{nn15}, one obtains
 a sufficiently
``singular" reaction term $q(f)$:
 \be
 \label{nn14}
 q(f) \sim (-y)^2 \sim f^{\frac 23} \asA y \to 0^-,
  \ee
i.e., $q(f)$ {\em is not} Lipschitz continuous at $f=0$, so that,
in the PDE setting, there occurs a typical problem of
non-uniqueness of solutions (if $q(f)$ is a source term for $f
\approx 0$).

However, we will show in Section \ref{SectLocR} that the
nonnegative TWs such as \ef{nn12} are not generic in the sense
that a.a. solutions are oscillatory and changing sign near finite
interfaces. In other words, nonnegativity of solutions is not an
invariant property of such quasilinear bi-harmonic flows in the
Cauchy problem setting (but can be in a FBP one, what was first
shown by Bernis--Friedman \cite{BF1} for quasilinear thin film
equations).
 Therefore, we
must neglect any opportunity to deal with nonnegative solutions of
such equations with arbitrary sufficiently smooth coefficients.

\ssk

For any $\l \ne 0$, the problem \ef{E5}, \ef{BC1} is of the
elliptic type, but {\em it is not variational}. Therefore, as in
\cite{GKPPI, GKPPIII},  one cannot use   advanced methods for
higher-order ODEs with potential operators associated with
homotopy-hodograph and other approaches \cite{PelTroy, KKVV00,
VV02} and/or Lusternik--Schnirel'man and fibering theory
\cite{GMPSobI, GMPSobII}. Therefore, the ODE \ef{E5}, though
looking rather simple, and, at least, simpler than most of related
fourth-order ODEs already studied in detail, represents a serious
challenge and cannot be tackled directly by known tools of modern
nonlinear analysis and operator theory.

\ssk

 Finally,
 we  study TWs for the {\em quasilinear tri-Harmonic equation}
\be
 \label{PME6}
  u_t=  (|u|^n u)_{xxxxxx}+ u(1-u) \inB \re \times \re_+ \whereA
  n>0.
   \ee


\subsection{Other related results of the present research}

The previous paper \cite{GKPPI} is devoted to the {\em semilinear
parabolic equations} for $n=0$, i.e.,
 \be
 \label{E4}
 u_t= -u_{xxxx} +u(1-u) \inB \re \times \re_+
  \LongA -\l f'=-f^{(4)}+f(1-f),
  \ee
and the boundary conditions now read
\be
 \label{BC1n0}
 f(y) \to 0   \andA f(y) \to 1 \asA y \to - \iy \quad
 \mbox{``maximally" exponentially fast.}
  \ee
This ``maximal" decay of $f(y)$ at infinity somehow includes some
kind of the remnants of a  ``minimality" of the possible TW
profiles, though any direct specification of such a property is
difficult to express rigorously.

\ssk

   Thus, the present paper and \cite{GKPPI} deal with higher-order {\em parabolic}
equations of reaction-diffusion type, while, in the third part of
this research \cite{GKPPIII}, we study KPP-type problems for
 other
types of PDEs including {\em dispersion}, {\em hyperbolic}, and
other ones.
 Namely,
in \cite{GKPPIII}, we will deal with higher-order hyperbolic and
dispersion equations such as
 \be
 \label{m3}
 u_{tt}= -u_{xxxx} +u(1-u) \andA  u_{ttt}= D_x^{(10)}u +u(1-u),
 \quad \mbox{etc.}
  \ee
As for more  ``exotic" PDE models, as a formal but quite
illustrative example, we consider, in \cite{GKPPIII}, higher-order
dispersion equations and end up with the following one:
 \be
 \label{m4}
 D_t^9 u= D_x^{11}u + u(1-u),
  \ee
  with eleventh-order ODE for the TW profiles
  \be
  \label{m41}
  -\l^9 f^{(9)}= f^{(11)} + f(1-f) \inB \re \quad (\mbox{plus\,\,
  (\ref{BC1})}).
  \ee
  We also consider in \cite{GKPPIII} an example of a {\em
  quasilinear} dispersion equation with a similar nonlinearity.


\subsection{Main questions to study}

 Thus, for quasilinear KPP-problems \ef{PME4} and \ef{PME6},
  the  main
 questions
  to study, here  and in \cite{GKPPI, GKPPIII} (though in the quasilinear cases
  the results achieved are weaker), are:

\ssk

 {\bf (I)} {\sc The problem of TW existence:} existence of travelling waves via
  analytical/numerical methods.
  Firstly, we easily
 show the positivity: $\l>0$ always.

\ssk


 {\bf (II)} {\sc The $\log t$-shift problem:} given some speed
 $\l_0>0$,
  derivation of a possible
 $\log t$-shifting of the moving front in the problem \ef{1.1},
 for classes of initial data $u_0(x)$ (for Heaviside data
 \ef{1.H}, we expect no shift and convergence to a TW; cf. some convincing stability results in \cite{Bert01, Bert99, Hu2008}
 for thin film equations). As in the semilinear cases, this phenomenon is also
  connected with a kind of an ``(affine) centre subspace
 behaviour" for the rescaled equation.


We thus show that  $\log t$-shifting phenomena, for some classes
of data) are also available for quasilinear degenerate KPP--4 and
other problems, though is less generic than for semilinear
higher-order equations, since require, in general, {\em still
unknown} hypotheses on $u_0(x)$.

\ssk

Unlike the semilinear analytic PDEs for $n=0$ in \cite{GKPPI}, the
rescaled quasilinear equations for $n>0$ are not analytic, so we
cannot proceed with a deeper study of the omega limits for such
equations.

\section{The basic higher-order KPP--4n problem}
\label{S2}

Consider the basic KPP--(4,1)$n$ (or simply KPP--4$n$, that cannot
confuse in the parabolic case) problem \ef{E4}, and let us begin
with its ODE counterpart \ef{E5}, \ef{BC1}.

\subsection{$\l$ is always positive}

 As in
\cite[\S~2]{GKPPI}, we first prove a  first simple result on the
positivity of admissible speeds $\l$:

\begin{proposition}
\label{Pr.Lambda1}
\be
\label{Lam1}
 \mbox{If there exists a solution $f(y) \not \equiv 0$
of $(\ref{E5})$, $(\ref{BC1})$, then $\l>0$.}
 \ee
 \end{proposition}

 \noi{\em Proof.} Multiplying the ODE \ef{E5} by $(|f|^n f)'$ and
 integrating by parts over $\re$  and using the boundary conditions \ef{BC1}
 yield
  \be
  \label{Lam2}
   \begin{aligned}
   \tex{
  -\l(n+1) \int |f|^n(f')^2}= &
 [F(f(y))]_{-\iy}^{+\iy} =- F(1)<0,\\
  &
  \tex{
  \mbox{where} \quad  F(f)= (n+1) \int_0^f |z|^nz(1-z)\,
 {\mathrm d}z. \qed
  }
  \end{aligned}
  \ee

\subsection{Numerical construction of TW profiles}

Again, as in \cite{GKPPI}, we begin with  presenting first
numerical results, which directly show the global structure of
such TW profiles to be, at least partially, justified
analytically. We use the {\tt bvp4c} solver of the {\tt MatLab}
with a sufficient accuracy; see more details in
\cite[\S~2]{GKPPI}. Note that, as the initial data for further
iterations, we always took the Heaviside function
 \be
 \label{Heav1}
 \mbox{initial data for numerical iterations are often} \quad
 H(-y),
  \ee
i.e.,  as in \ef{1.H}. This once more had to help us to converge
to a proper ``minimal" profile (indeed, there are many other TW
profiles), though, of course ,this was not guaranteed {\em a
priori}. We kept this rule for all other KPP--$(k,l)$ problems of
interest in \cite{GKPPI, GKPPIII}.

For convenience, we perform all the calculations for the function
 \be
 \label{FF11}
  \tex{
 F(y)=|f(y)|^n f(y) \LongA
   F^{(4)}= \frac \l{n+1}\, |F|^{-\frac n{n+1}}
  F'+|F|^{-\frac n{n+1}}F \big(1-|F|^{-\frac n{n+1}}F\big),
 }
  \ee
so that the ODE in \ef{E5} becomes {\em semilinear}.

 \ssk

\begin{figure}
\centering
\includegraphics[scale=0.85]{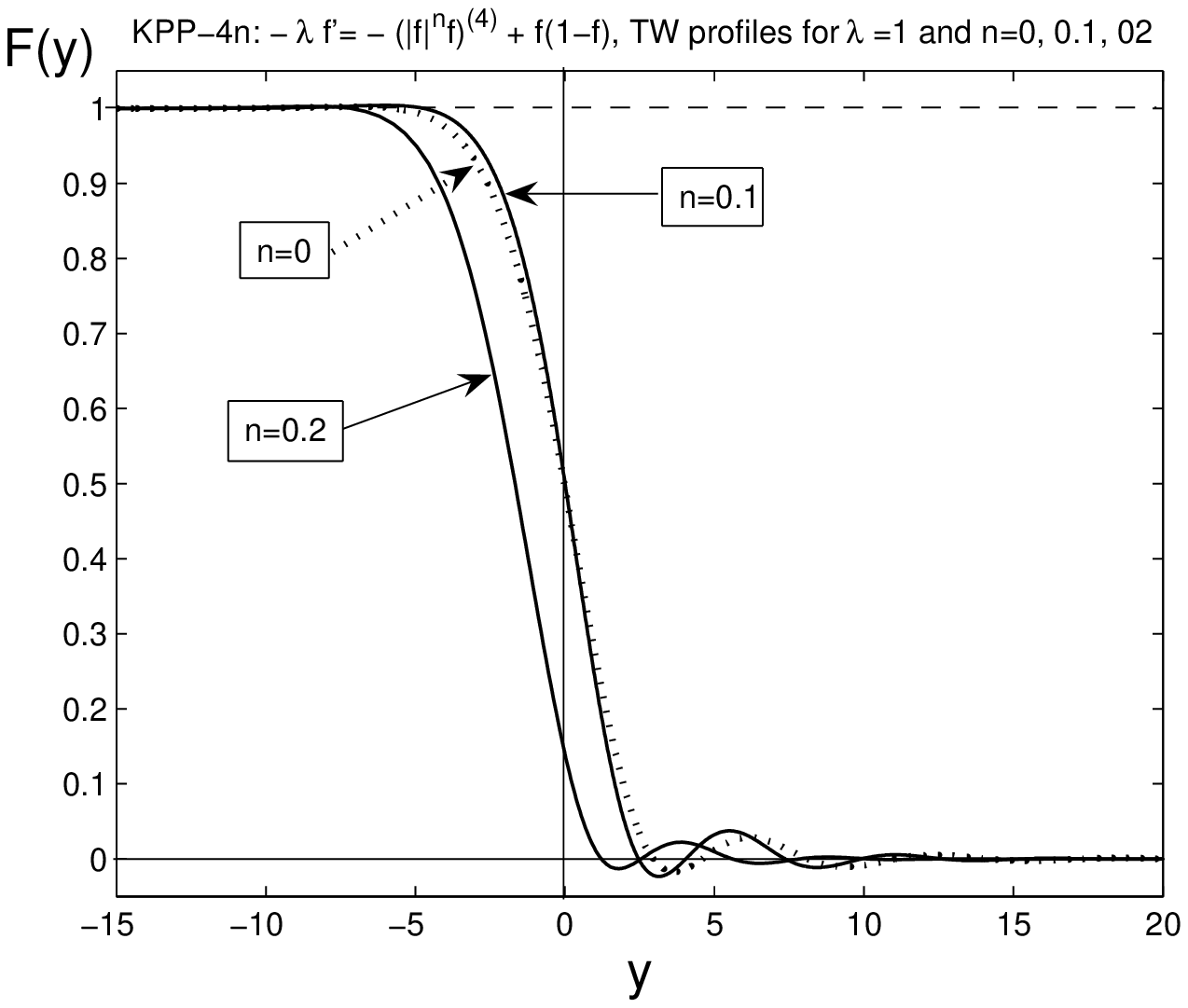}  
\vskip -.3cm
  \caption{The TW profiles $F(y)$ of
(\ref{FF11}), \ef{BC1} for
 $\l=1$, $n=0, 0.1, 0.2$.}
 \label{Fn1}
\end{figure}

\begin{figure}
\centering
\includegraphics[scale=0.85]{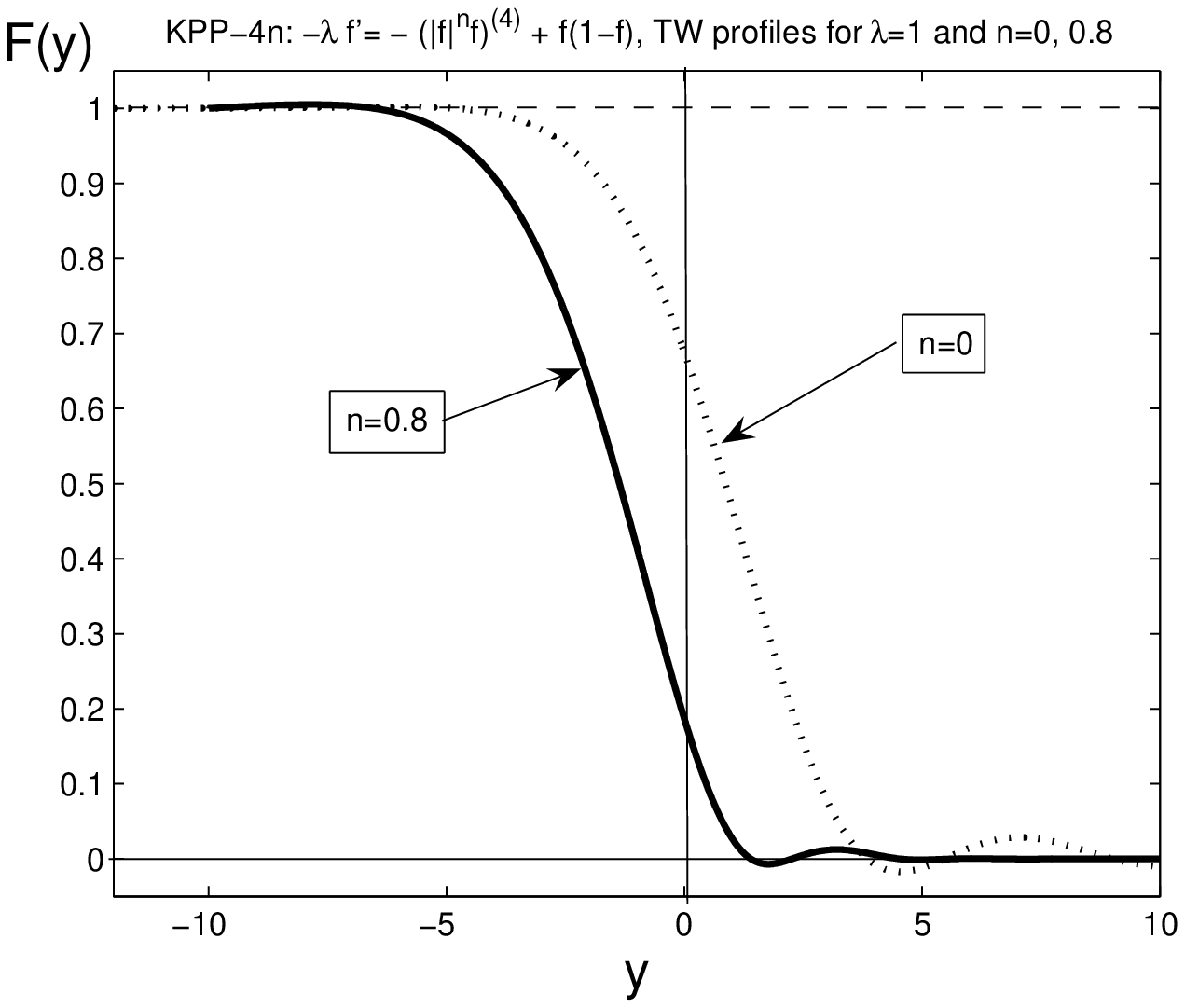}  
\vskip -.3cm
  \caption{The TW profiles $F(y)$ of
(\ref{FF11}), \ef{BC1} for
 $\l=1$, $n=0.8$.}
 \label{Fn2}
\end{figure}

Figure \ref{Fn1} shows TW profiles $F$ for
 $\l=1$
for
 sufficiently small $n=0, 0.1, 0.2$. The semilinear case $n=0$ was
 used to get a good comparison with the results in
 \cite[\S~2]{GKPPI}.
  For larger $n=0.8$, with the same $\l=1$, the  numerical results  are shown in Figure
  \ref{Fn2}.


\begin{figure}
\centering
\includegraphics[scale=0.85]{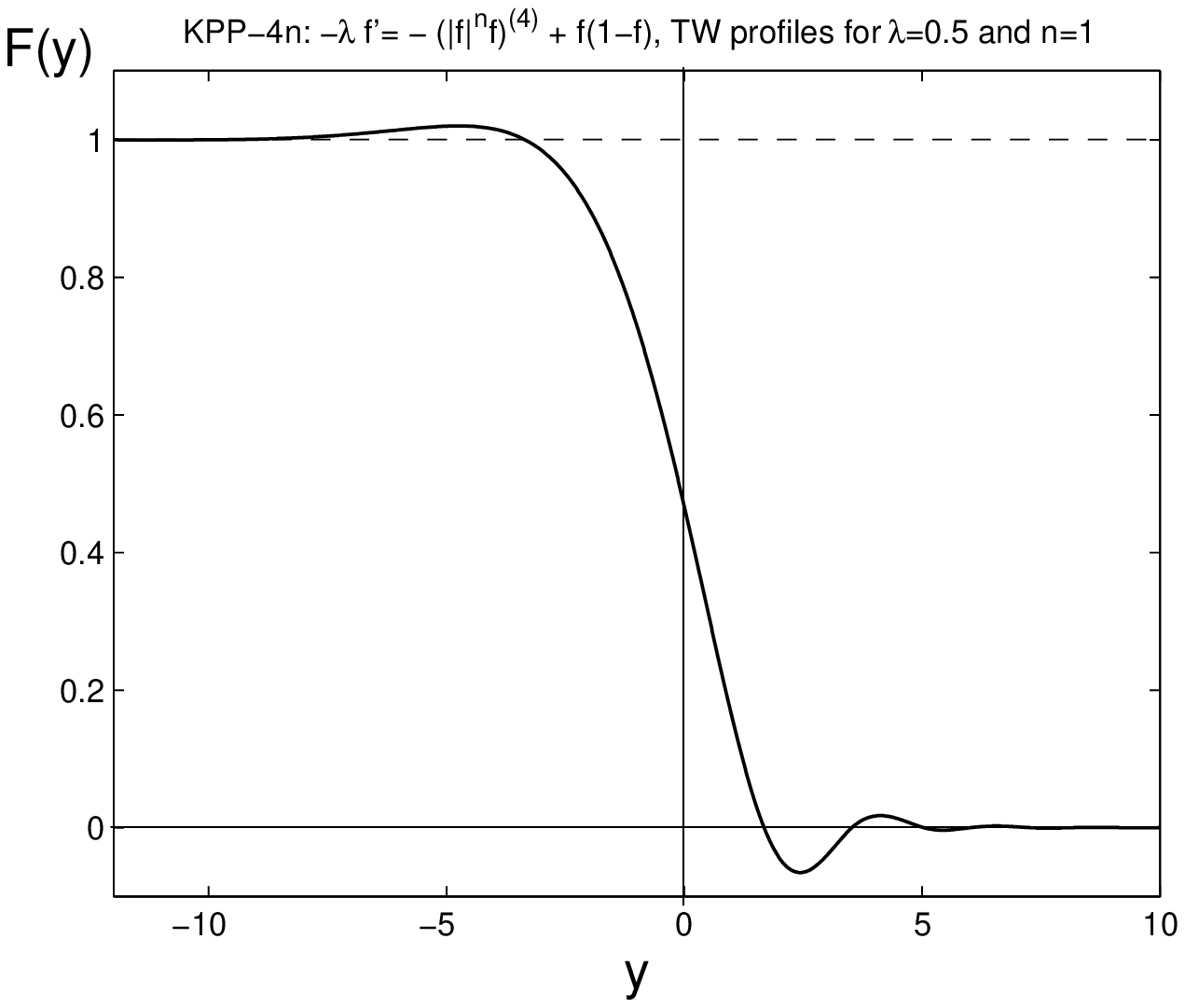}  
\vskip -.3cm
  \caption{The TW profiles $F(y)$ of
(\ref{FF11}), \ef{BC1} for
 $\l=0.5$, $n=1$.}
 \label{Fn3}
\end{figure}

\begin{figure}
\centering \subfigure[oscillations about $F=0$]{
\includegraphics[scale=0.52]{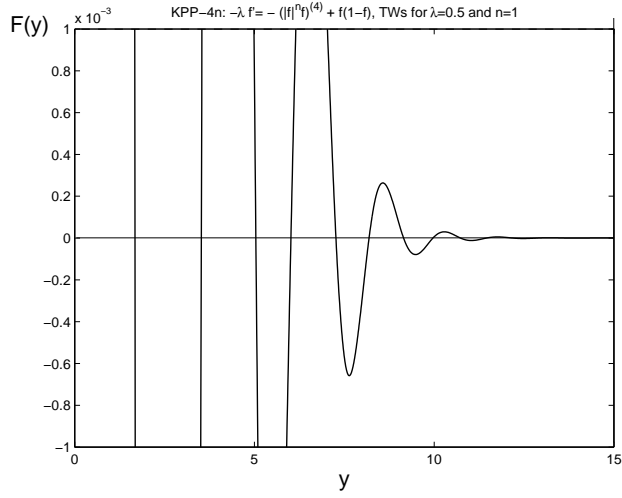}             
} \subfigure[oscillations about $F=1$]{
\includegraphics[scale=0.52]{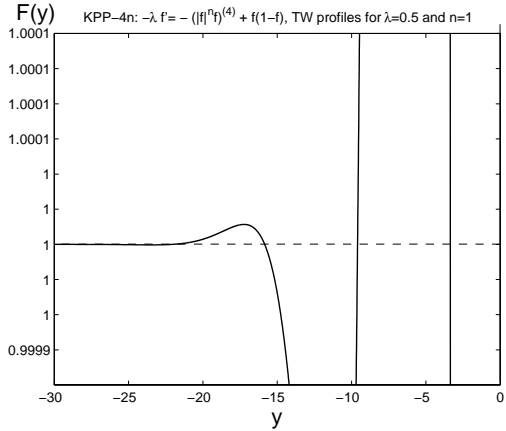}                        
}
 \vskip -.3cm
\caption{\rm\small An oscillatory
 convergence of the TW profile $F(y)$ to 0, $y \gg1$ (a) and to 1, $y \ll
-1$; $\l=0.5$ and $n=1$.}
 \label{Fn4}
\end{figure}

\begin{figure}
\centering
\includegraphics[scale=0.85]{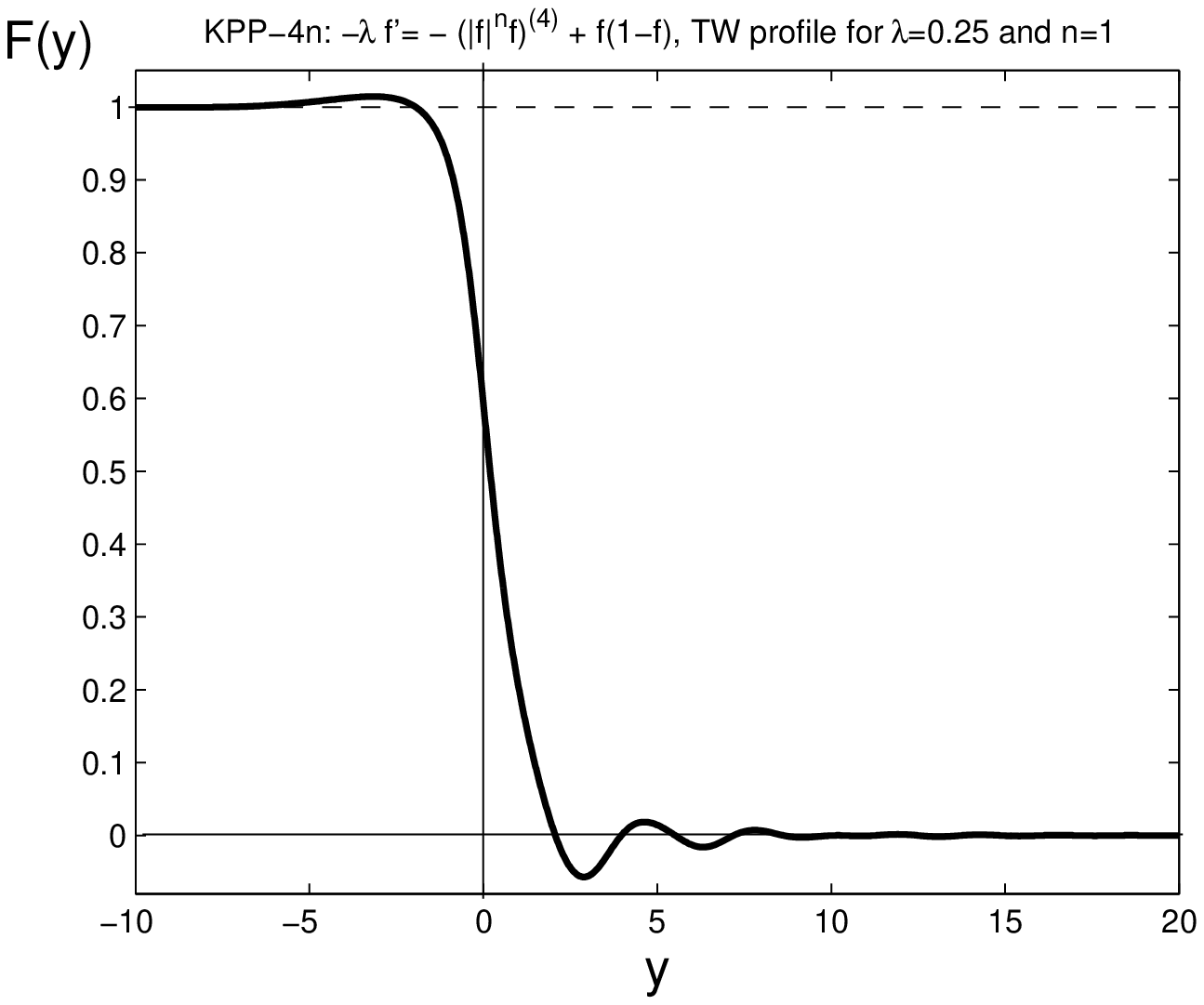}  
\vskip -.3cm
  \caption{The TW profiles $F(y)$ of
(\ref{FF11}), \ef{BC1} for
 $\l=0.25$, $n=1$.}
 \label{Fn5}
\end{figure}

We next fix $n=1$ and $\l=0.5$ as Figure \ref{Fn3} shows. A more
detailed structure of oscillations of solutions about equilibria
$F=0$ (a) and $F=1$ (b) is presented in Figure \ref{Fn4}.
Separately, Figure \ref{Fn5} shows TW profiles $F(y)$ as in
\ef{FF11} for $n=1$ and $\l= \frac 14$.


Finally, we must admit that, for $n>0$, and, especially, for
larger values $n \sim 1$, it is much more difficult to get
reliable well-converging numerical results unlike the more
straightforward semilinear case $n=0$, \cite{GKPPI}. Many our
attempts, even rather time-consuming (lasting for several hours of
the {\tt MatLab}) led to oscillatory non-converging profiles
and/or to singular Jacobians (a kind of ``blow-up" of
computations).

However, our numerical experiments produced a rather convincing
evidence of existence of TW profiles for a variety of $n>0$ and
velocities $\l>0$, which we  will need to develop further.

Thus, we begin with necessary local analysis of the behaviour of
TW profiles $f(y)$ close to equilibria $f=1$ and $f=0$ (see the
next section), in order to fix  dimensions of their stable and
unstable manifolds, for further matching.

\section{The 2D stable bundle as $y \to -\iy$ and instabilities}
 \label{Sy1}

This analysis is not much different from that for $n=0$ in
\cite[\S~2.3]{GKPPI}.

 \subsection{Linearization about $f=1$}

 Thus, setting $f=1+g$ in \ef{E5} and linearizing
yield the following characteristic equation:
\be
\label{t9}
 - \l g'= (n+1)g^{(4)} + g \andA
  g(y)=\eee^{\mu y} \LongA H_-(\mu,\l) \equiv (n+1)\mu^4 -\l \mu +1 =0.
   \ee
Therefore, for $\l=0$, we have 2D stable (${\rm Re}\,(\cdot) >0$)
and unstable (${\rm Re}\,(\cdot) <0$) bundles with the roots
 \be
 \label{t10}
  \tex{
   \mu_\pm(0)= \pm \frac {1+\ii}{\sqrt 2}(n+1)^{-\frac 14} \andA
\bar\mu_\pm(0)= \pm \frac {1-\ii}{\sqrt 2} (n+1)^{-\frac 14}. }
 \ee
 By continuity, for all small $\l>0$, there exists 2D stable
 manifold of the equilibrium 1 with the roots
  \be
  \label{root11}
   \tex{
 \mbox{stable bundle}: \quad  \mu_+(\l) \andA \bar \mu_+(\l).
  }
  \ee

\begin{proposition}
 \label{Pr.osc2} For any $n \ge 0$,
 at least for all small $|\l|>0$, the linearized equation
 $\ef{t9}$
$($and hence the KPP--$4n$ one in $\ef{E5}$ for $f \approx 1$ for
$y \ll -1)$ admits a $2D$ stable family of oscillatory solutions
as $y \to -\iy$, and a $2D$ unstable one of exponentially
divergent orbits.


 \end{proposition}



\subsection {Local blow-up and other instabilities to $-\iy$}

In order to verify the global continuation properties of stable
bundles, one need to check whether the {\em nonlinear} ODE \ef{E5}
admits blow-up and the dimension of such an unstable manifold. To
this end, we re-write it down and, as usual, neglect the linear
lower-order terms, which by standard local interior regularity are
negligible for $|f| \gg1$:
 \be
 \label{t12}
  (|f|^n f)^{(4)}=-f^2 + f +\l f' = - f^2(1+o(1)) \asA f \to \iy.
   \ee
For any $n \in [0,1)$, the unperturbed equation has the following
exact blow-up solution: for the function $F$ in \ef{FF11},
 \be
\label{t13}
 \tex{
  F^{(4)}=-|F|^{\frac 2{n+1}} \LongA  F_0(y)=-
  C_0(n)(Y_0-y)^{-\frac {4(n+1)}{1-n}}
\to - \iy \asA y \to Y_0^-,
 }
  \ee
  where $Y_0 \in \re$ is a fixed arbitrary blow-up point and $C_0(n)>0$ is an easily computed constant.
   Studying the
  dimension of its stable manifold, as in \cite[\S~2.4]{GKPPI},
  yields


 \begin{proposition}
 \label{Pr.Blow1}
 For the ODE $\ef{E5}$, with any $n \in [0,1)$, the stable manifold of blow-up solutions is
 1D depending on a single parameter being their blow-up point $Y_0 \in
 \re$.
  \end{proposition}

For $n=1$, the asymptotic ODE \ef{t13} has the form
 \be
 \label{t131}
 F^{(4)}=-|F| \LongA F_0(y)=- \eee^{-y} \to - \iy \asA y \to -\iy,
  \ee
  so, instead of having finite $y$-blow-up, one gets this
 negative  exponential growth for $y \ll -1$, with a similar
 one dimension.

 Finally, for $n>1$, we have algebraically growing solutions:
 \be
 \label{t132}
 F_0(y)= -C_1 (-y)^{\frac {4(n+1)}{n-1}} \to -\iy \asA y \to -
 \iy \quad (C_1>0).
  \ee

\section{Oscillatory solutions
 near finite interfaces: 3D asymptotic bundle}
  \label{SectLocR}

  This part of the asymptotic analysis is essentially different
  from that for the semilinear case $n=0$ in \cite[\S~2]{GKPPI};
  however, such results are already known \cite{GalRDE4n}, so we
  can omit some details.

  Thus, for $n>0$, we describe the generic oscillatory behaviour
  of solutions of (\ref{FF11}) close to finite interfaces.
  According to pioneering results in \cite{Bern88, BMc91}, since
  1988, it was known that
  ODEs like (\ref{FF11}) admit compactly supported solutions
 of changing sign near finite interfaces. Let $f(y)$ vanish at the interface $y
 \to y_0>0$, so that $f(y) \equiv 0$ for $y > y_0$.
 Then, making for convenience the reflection $y \mapsto y_0-y$,
 with $y>0$ small enough, and keeping  the
 leading first two terms in (\ref{FF11}) for $y \approx 0^+$, after integration once
 we obtain an exponentially small perturbation of the following
 third-order equation (the sign ``$-$" on the right-hand side appears because of the reflection in $y$):
  \be
  \label{FF98}
  F'''= - \l |F|^{- \frac n{n+1}}F, \quad y>0, \quad F(0)=0.
  \ee
We next scale out the positive constant $\l$ to get the ODE
  \be
  \label{le2N}
   F'''=- \big|F\big|^{-\frac
  n{n+1}}F, \quad y>0.
   \ee

 We next  describe
oscillatory solution of changing sign of the ODE (\ref{le2N}),
with zeros concentrating at the given interface point $y=0^+$. Let
us mention again that oscillatory properties of solutions are a
common feature of higher-order degenerate ODEs. We refer to first
results in \cite{Bern88, BMc91, {GalRDE4n}}, to \cite{Gl4, GBl6}
(thin film equations), and to \cite[Ch.~3-5]{GSVR}, where further
examples can be found.

 To this end,
by the scaling invariance of (\ref{le2N}), we look for its
solutions of the form
 \be
 \label{le3}
  \mbox{$
 F(y) = y^\mu \varphi(s), \quad s= \ln y, \quad
 \mbox{where} \,\,\, \mu = \frac {3(n+1)}{n}> 3 \,\,\,\mbox{for} \,\,\, n>0,
  $}
  \ee
  where $\varphi(s)$ is  the so-called
  {\em oscillatory component}.
 Substituting (\ref{le3}) into (\ref{le2N}) yields
 the following third-order equation for $\varphi(s)$:
 \be
 \label{le4}
  P_3(\var)= - |\var |^{-\frac
  n{n+1}} \var,
  \ee
   where $P_k$ denote linear differential polynomials
   obtained by a simple recursion procedure (see
   \cite[p.~140]{GSVR}), so that
 $$
 \begin{matrix}
P_1(\varphi)=\var'+ \mu \var, \quad P_2(\var)= \var''+(2\mu-1)
\var'+ \mu(\mu-1)\var, \smallskip\smallskip \\
  P_3(\varphi)=\var'''+ 3(\mu-1) \var'' + (3 \mu^2- 6 \mu +2) \var'
 + \mu(\mu-1)(\mu-2) \var.
 \end{matrix}
 $$

According to (\ref{le3}), we are interested in uniformly bounded
global solutions $\varphi(s)$ that are well defined as $s= \ln y
\to -\infty$, i.e., as $y \to 0^+$. The best candidates for such
global orbits of (\ref{le4}) are periodic solutions $\varphi_*(s)$
that are defined for all $s \in \re$. Indeed, they can describe
suitable (and, possibly, generic) connections with the interface
at $s=-\infty$. See \cite{GalRDE4n} for the following and other
related results.

\begin{proposition}
 \label{Pr.Per}

For any $n>0$,  $(\ref{le4})$
 has a periodic solution $\varphi_*(s)$ of changing sign.


 \end{proposition}


Two problems remain open:

(i) uniqueness of the periodic solution $\var_*(s)$, and

(ii) stability (hyperbolicity) of $\var_*(s)$ as $s \to + \infty$.

\noi Numerically, we have obtained  positive answers to both
questions. In particular, (i) and (ii) imply that there exists a
unique (up to translation) periodic bounded connection with $s=
-\infty$, where the interface is situated.

\ssk

\begin{figure}
\centering \subfigure[$n=0.75$]{
\includegraphics[scale=0.52]{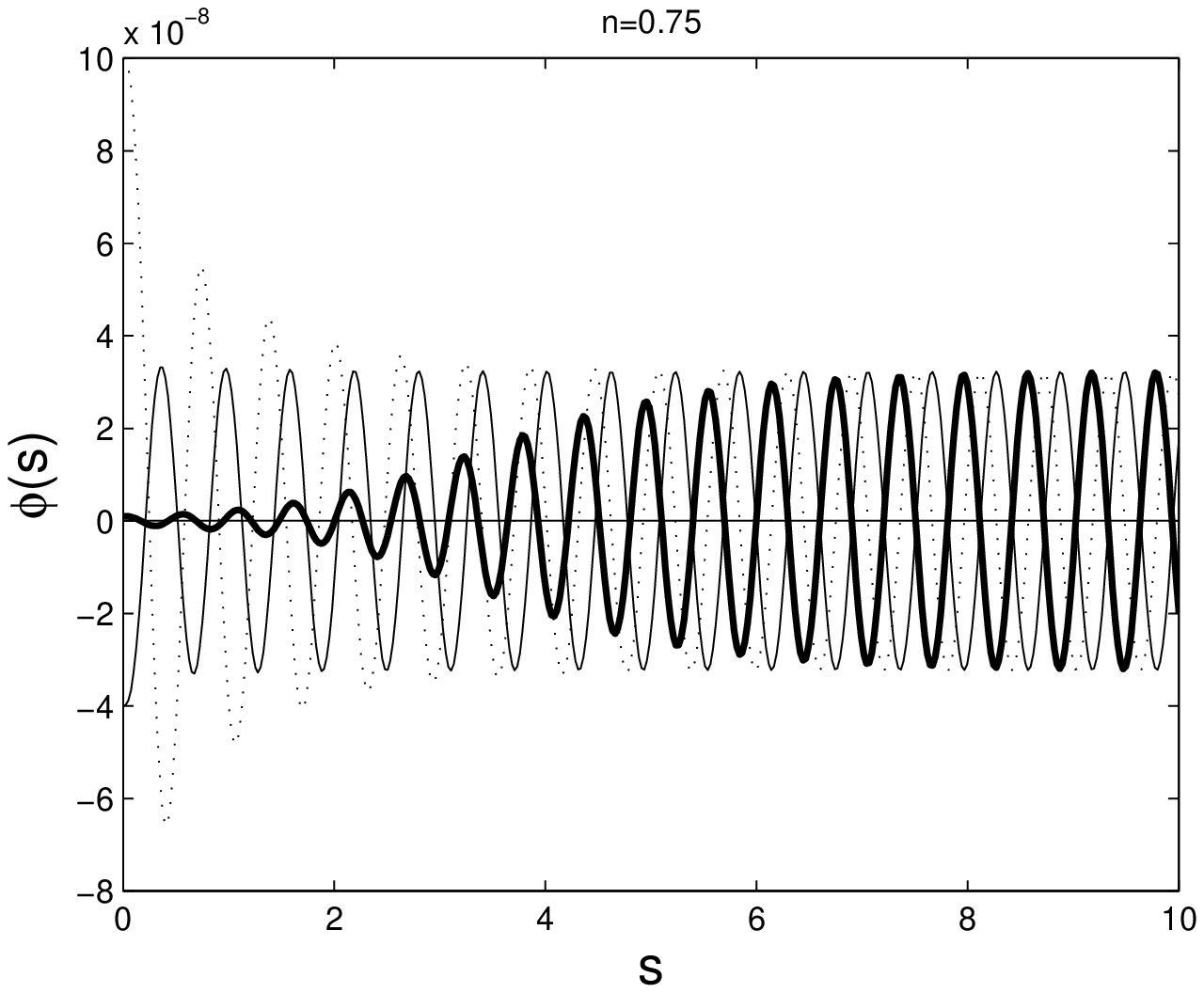}
} \subfigure[$n=1$]{
\includegraphics[scale=0.52]{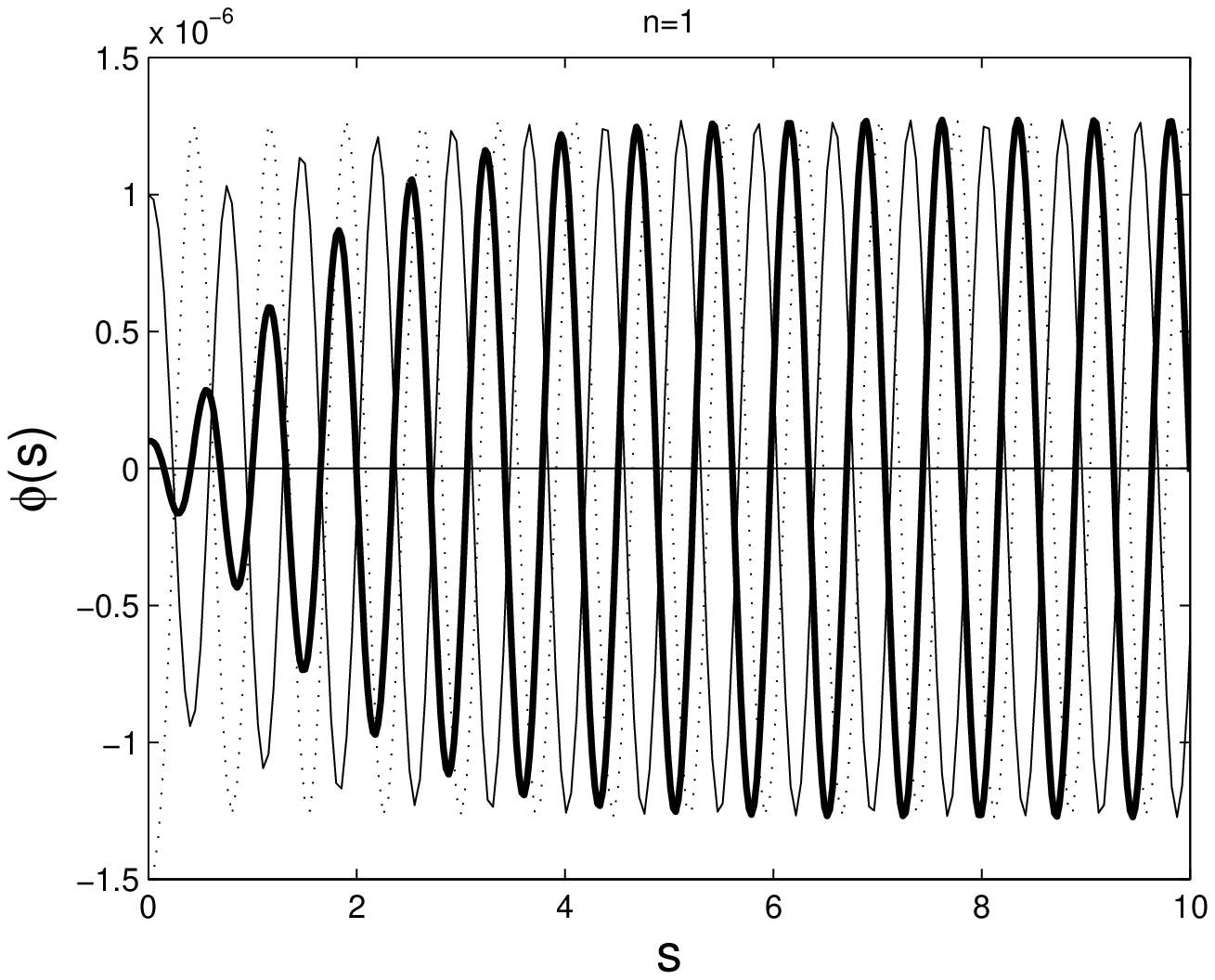}
}
 \subfigure[$n=2$]{
\includegraphics[scale=0.52]{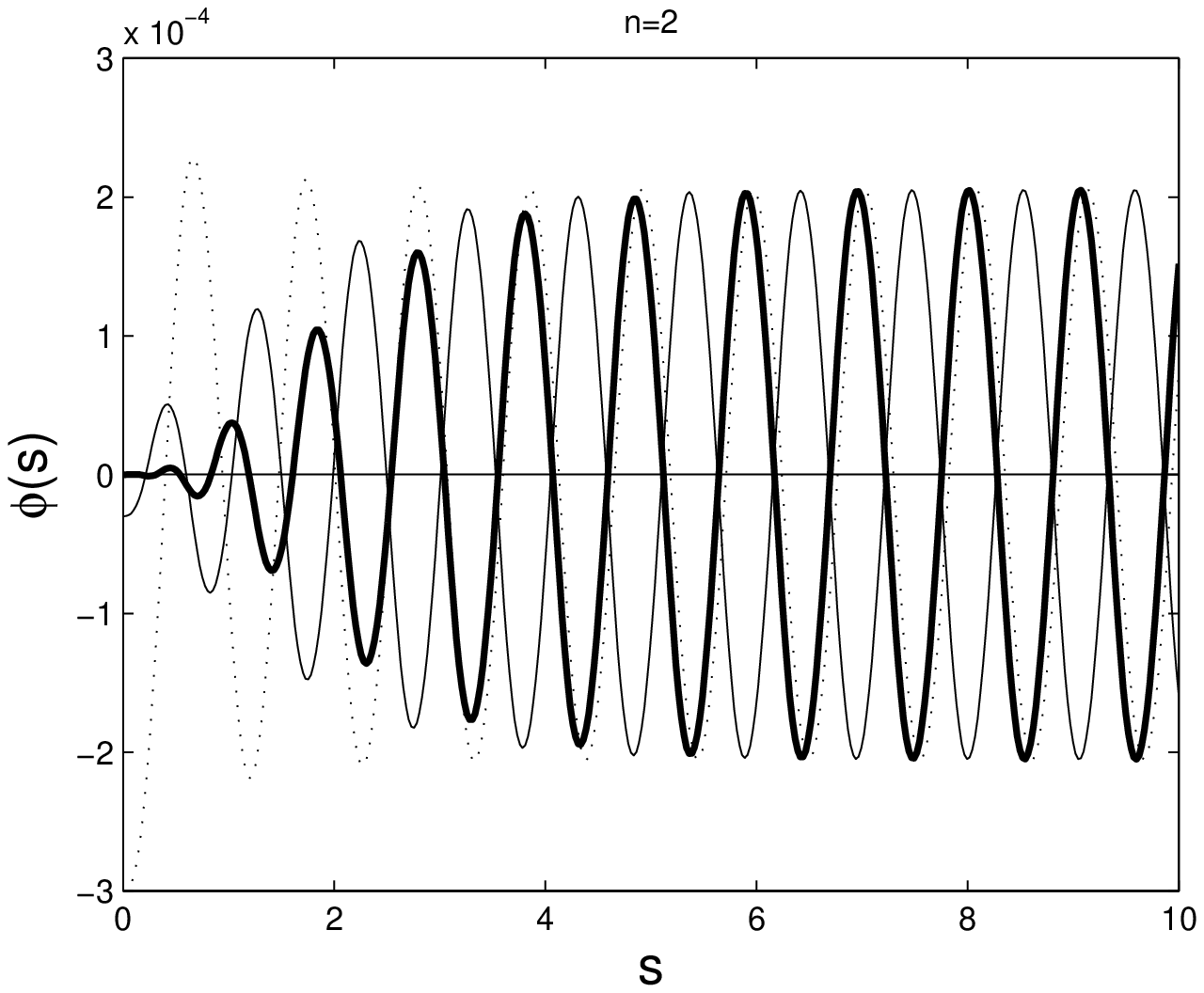}
} \subfigure[Large $n$]{
\includegraphics[scale=0.52]{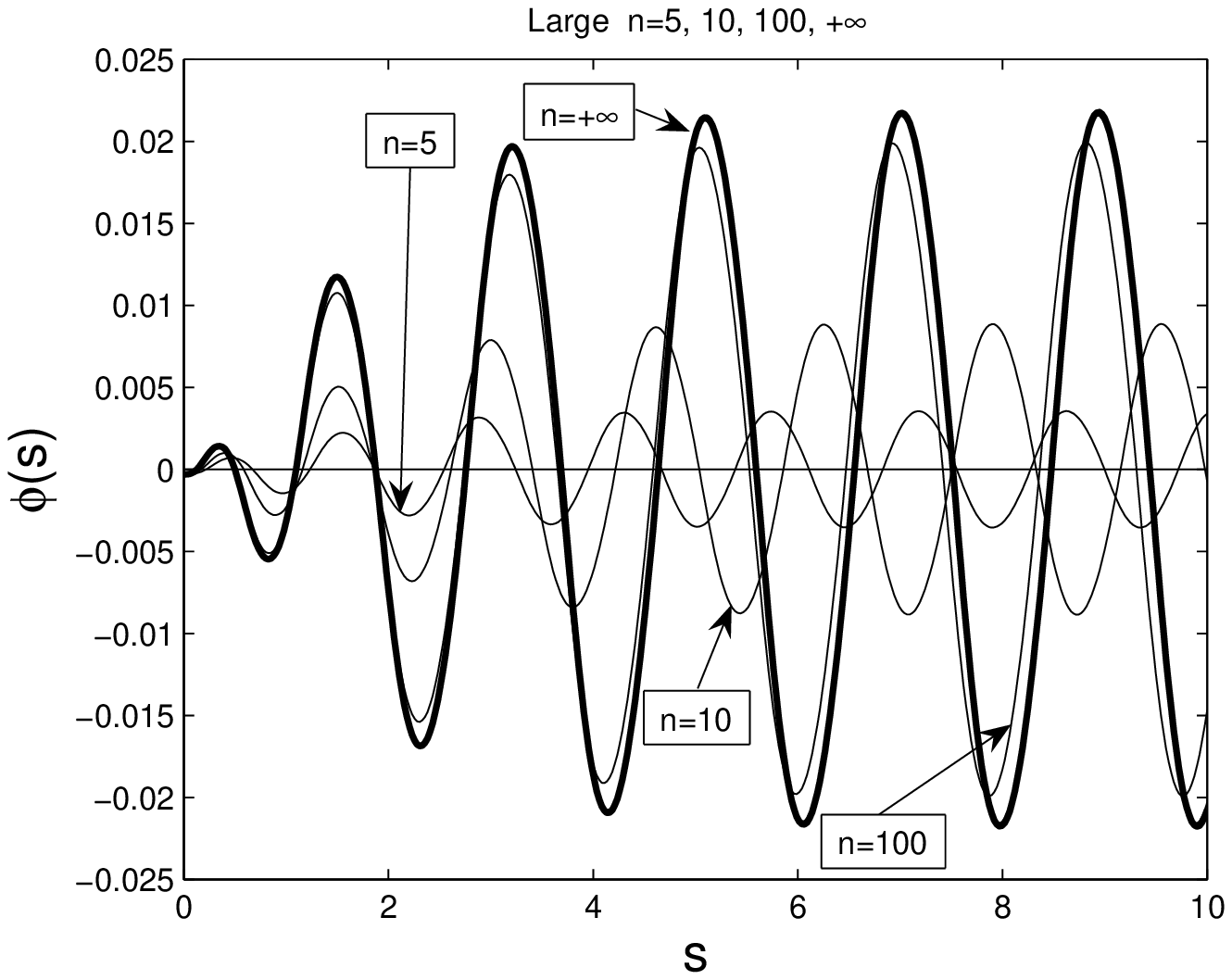}
}
 \vskip -.4cm
\caption{\rm\small Convergence to a stable periodic orbit of the
ODE (\ref{le2N}) for various $n>0$.}
 \label{FOsc2}
\end{figure}

 The convergence to the unique stable periodic
 solution of (\ref{le4}) is shown in Figure \ref{FOsc2} for
 various $n > 0$.
 Different  curves therein correspond to different Cauchy data
$\varphi(0)$, $\varphi'(0)$,  $\varphi''(0)$  prescribed at $s=0$.
For $n$ smaller than $\frac 34$, the oscillatory component gets
extremely small, so an extra scaling is necessary, which is
explained in \cite[\S~7.3]{Gl4}. A more accurate passage to the
limit $n \to 0$ in (\ref{le4}) is done there in Section 7.6 and in
Appendix B. This  explains the continuous deformation as $n \to 0$
of oscillatory structures in (\ref{le3}) to linear ones in the
exponential tail of the  kernel $F(y)$ of the fundamental solution
of the bi-harmonic operator $D_t+D_x^4$.

In (d), we also present  the periodic solution for $n=+\infty$,
where (\ref{le4}) takes a simpler form (see an algebraic
construction of the unique periodic solution in
\cite[\S~7.4]{Gl4})
 $$
 P_3(\varphi)=- {\rm sign} \, \varphi.
 $$


Finally, given the 
periodic $\varphi_*(s)$  of (\ref{le4}),
 as a natural
way to approach the interface point $y_0=0 $ according to
(\ref{le3}), we have that the ODE (\ref{le2N})
  generates at the singularity set $\{f=0\}$
 \be
 \label{as55}
\mbox{a 3D local asymptotic bundle with parameters $y_0$, a phase
shift in $s \mapsto s + s_0$},
 \ee
  and the parameter $\e>0$ of the scaling group for the ODE
  \ef{le2N}.
 \be
 \label{sc1}
 F \mapsto \e^{\frac{3(n+1)}n} \,F, \quad y \mapsto \e \, y \quad
 (\e>0).
  \ee
  Notice that this scaling invariance has been lost in the
  approximate ODE (\ref{le2N}).

\ssk


\noi{\bf Remark: nonnegative solutions for $\l<0$ (Example 3).}
One can see that \ef{le4} does not admit nonnegative solutions and
equilibria. However, if $\l<0$ in \ef{FF98}, then scaling out such
$\l$ leads to ODEs \ef{le2N} and \ef{le4} with the opposite sign
on the right-hand side, so that there exists a positive
equilibrium
 \be
 \label{pos1}
  \tex{
 P_3(\var)=|\var|^{- \frac n{n+1}} \var
 \LongA \var_0=\big[ \mu(\mu-1)(\mu-2)]^{-\frac{n+1}n}
 }
 \ee
(the rest of solutions are still oscillatory). Existence of such a
particular positive solution in \ef{pos1} is the ``local" origin
of the global nonnegative TW in Example 3 in Section \ref{Sect1},
though such a behaviour is not structurally stable for both ODEs
involved. Note also that the calculation \ef{Lam2} does not apply
to \ef{nn11} with a more complicated reaction $q(f)$.


\subsection{Well-defined matching of stable manifolds of equilibria 0 and 1}

Similar to \cite{GKPPI} for $n=0$ (when the conclusion is indeed
more straightforward), we thus observe that, for $n>0$ and
sufficiently small $\l>0$, there is a well defined matching
procedure of stable manifolds of equilibria $f=0$ (see \ef{as55})
and $f=1$ (Proposition \ref{Pr.osc2}). The overall relation of the
dimensions is
 \be
 \label{dim1}
 3_{(y \gg 1, \, f \to 0)} - 1 =2_{(y \ll -1, \, f \to 1)},
 \ee
 where ``$-1$" stands for the dimension of the unstable (blow-up for $n<1$)
manifold in Proposition \ref{Pr.Blow1}. In other words, \ef{dim1}
implies a well posed algebraic system of two equations with two
unknowns. In the case of analytic manifolds as in \cite{GKPPI},
this guarantees at most a countable number of solutions (or a
finite one of uniformly bounded TW profiles).

In the present degenerate case $n>0$, any analytic dependence on
parameters is also plausible but difficult to prove. However, we
think that \ef{dim1} confirms that, for small $\l>0$, the family
of TW profiles is always discrete (what have been seen in
numerical experiments), and then such (most probably, finite)
number of $\l$-branches can be extended for larger $\l>0$ by
classic nonlinear integral operator theory \cite{KrasZ, Deim}.
Note that a fully global extension in $\l$ is hardly possible: as
was shown in \cite{GKPPI} for $n=0$, there is a clear phenomenon
of existence of a maximal speed $\l_{\rm max}$ such that, for $\l
> \l_{\rm max}$, TW profiles $f(y)$ are nonexistent. By a
continuity argument, we expect that a similar phenomenon takes
place for $n>0$, but this is difficult to check even numerically.

\section{A few words on ``more quasilinear" KPP--4$n$ problem}
 \label{S2Quas}

 This is about the KPP-setting for the following quasilinear
 equation:
 \be
  \label{qu1}
  u_t=-(|u|^n u)_{xxxx} + |u|^n u(1-|u|^n u) \LongA
  -\l f'=-(|f|^n f)^{(4)}+|f|^nf(1-|f|^nf),
   \ee
where the source  contains the same nonlinear term $|u|^n u$ as
the porous medium diffusion one. A one advantage of this model is
that the corresponding ODE for $F=|f|^n f$ is less complicated
than that in \ef{FF11}:
 \be
\label{qu2}
 \tex{
  F^{(4)}= \frac \l{n+1}\, |F|^{-\frac n{n+1}} F' +
F(1-F), }
 \ee
 which, as might be expected, could improve convergence of
 numerical schemes. However, this did not happen, and convergence,
 though being slightly better, was not improved somehow
 essentially.

\begin{figure}
\centering
\includegraphics[scale=0.85]{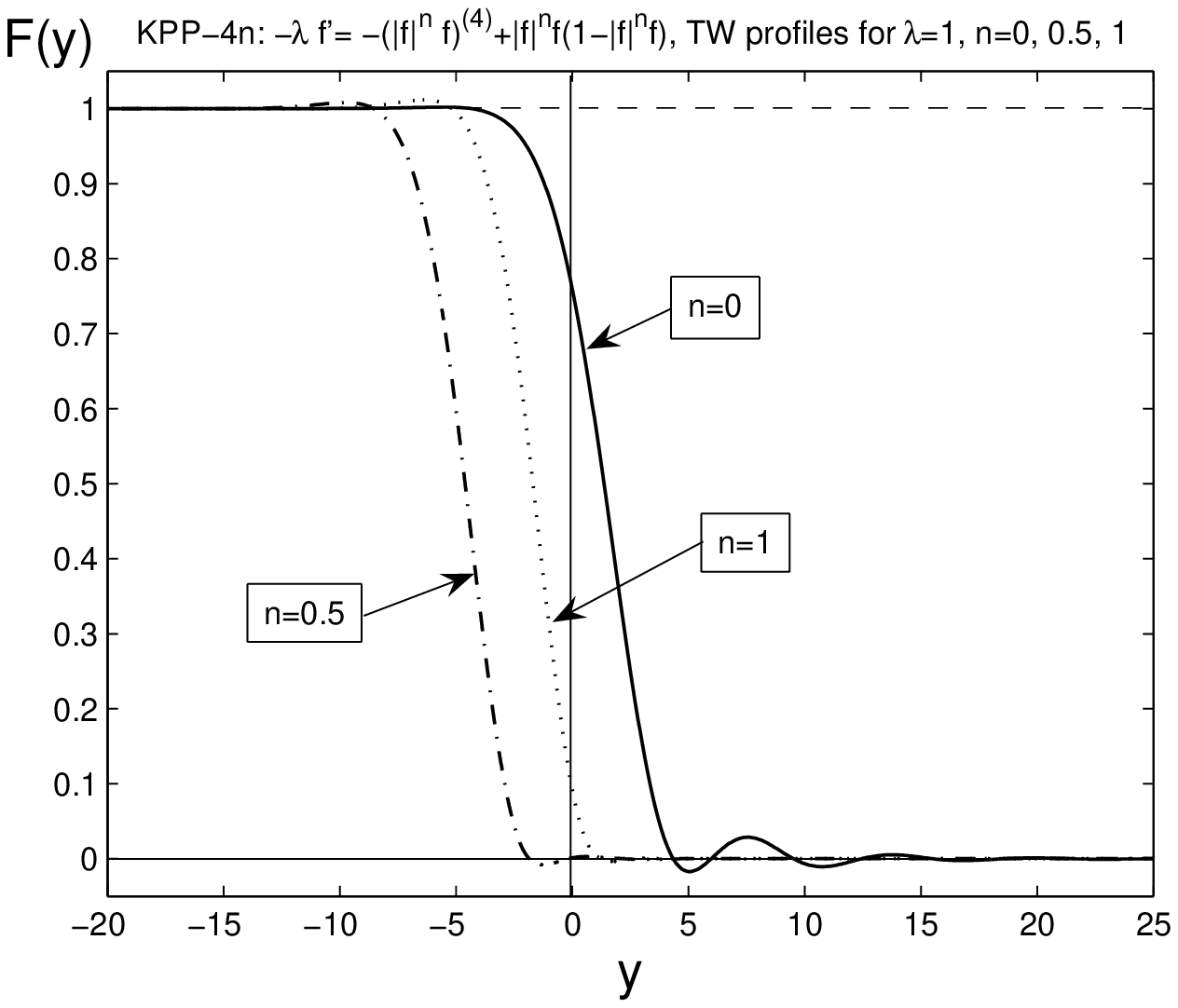}  
\vskip -.3cm
  \caption{The TW profiles $F(y)$ of
(\ref{qu2}), \ef{BC1} for
 $\l=1$, $n=0, 0.5$ and $1$.}
 \label{FQ1}
\end{figure}

\begin{figure}
\centering
\includegraphics[scale=0.85]{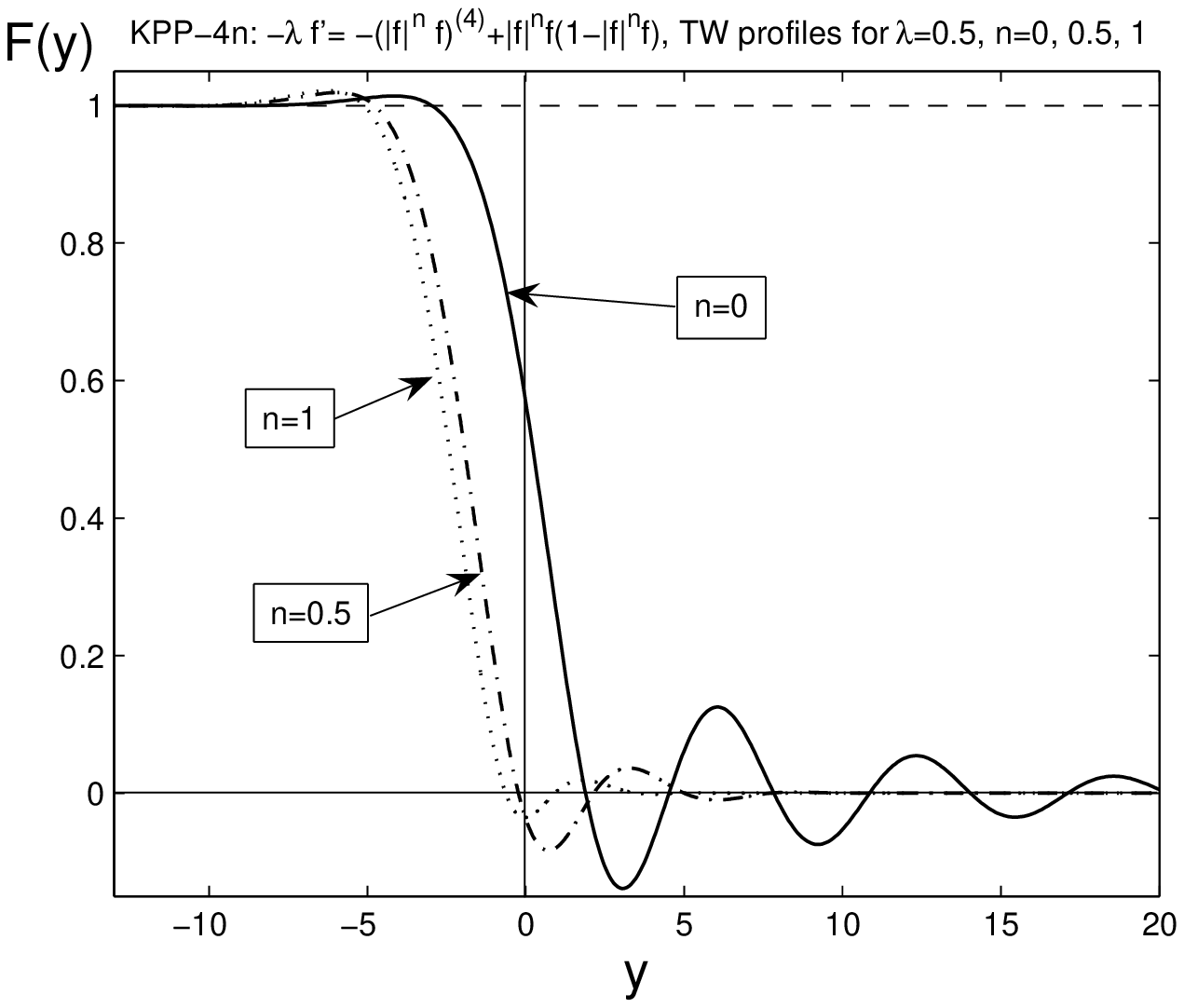}  
\vskip -.3cm
  \caption{The TW profiles $F(y)$ of
(\ref{qu2}), \ef{BC1} for
 $\l=0.5$, $n=0, 0.5$, and $1$.}
 \label{FQ2}
\end{figure}

 In Figures \ref{FQ1} and \ref{FQ2}, we show TW profiles $F(y)$
 for $n=0, \, 0.5$, and 1 for the two cases $\l=1$ and $\l=0.5$.


\begin{figure}
\centering
\includegraphics[scale=0.85]{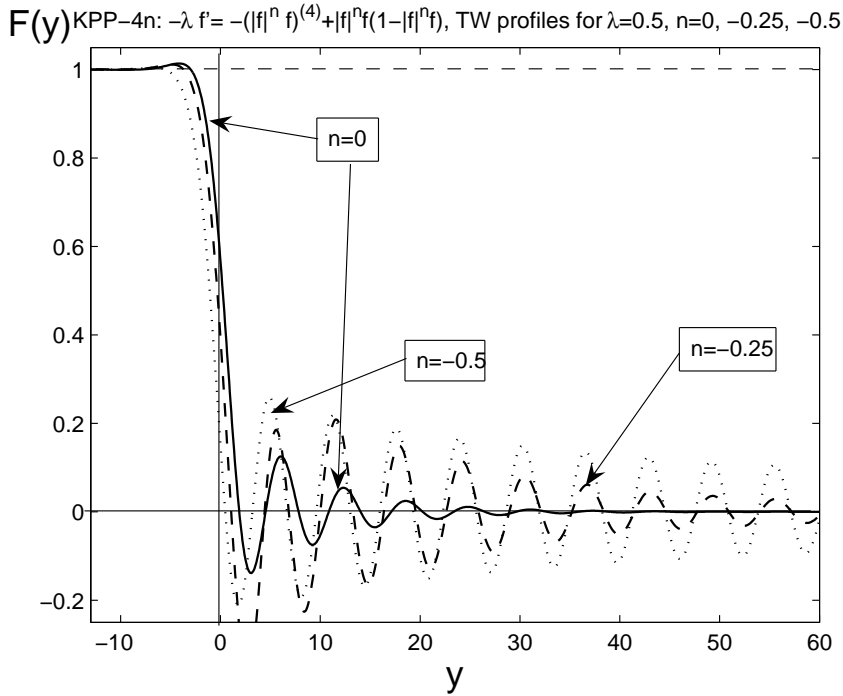}  
\vskip -.3cm
  \caption{The TW profiles $F(y)$ of
(\ref{qu2}), \ef{BC1} in the ``fast diffusion" case for
 $\l=0.5$, $n=0, -0.25$ and $-0.5$.}
 \label{FQ3}
\end{figure}

An unusual application of such a quasilinear bi-harmonic equation
\ef{qu1} is that it allows to treat the ``fast diffusion" case
$n<0$; more precisely, $n \in (-1,0$). In Figure \ref{FQ3}, we
show such oscillatory profiles for $n=0$, $-0.25$, and $-0.5$ for
$\l=0.5$. The convergence for $n<0$ is much better than in the
standard porous medium case $n>0$. It is worth mentioning that,
since for $u \approx 0$ the source term becomes a non-Lipschitz
function,
 $Q(u) \sim +|u|^n u$, $n<0$, there appears a standard
 non-uniqueness of solutions even in the simple ODE:
  \be
  \label{qu5}
  u_t=|u|^n u, \quad t>0, \quad u(0)=0 \LongA \exists \,\,\, u(t)=(|n|t)^{\frac
  1{|n|}}>0,
   \ee
   together with the trivial solution $u(t) \equiv 0$. In such
   cases, in PME and fast diffusion theory with blow-up, extinction, quenching, and other singularities, an extended semigroup
   of proper extremal (maximal or minimal) solutions is
   constructed; see \cite[Ch.~6,7]{GalGeom} and references therein.
   However, such a theory is fully nonexistent for higher-order
   parabolic flows with such singularities. Creating such a theory
   for equations with no Maximum Principle is a hard open problem.


\section{A few words on non-divergent KPP thin film problem}
 \label{STFE}

\begin{figure}
\centering
\includegraphics[scale=0.85]{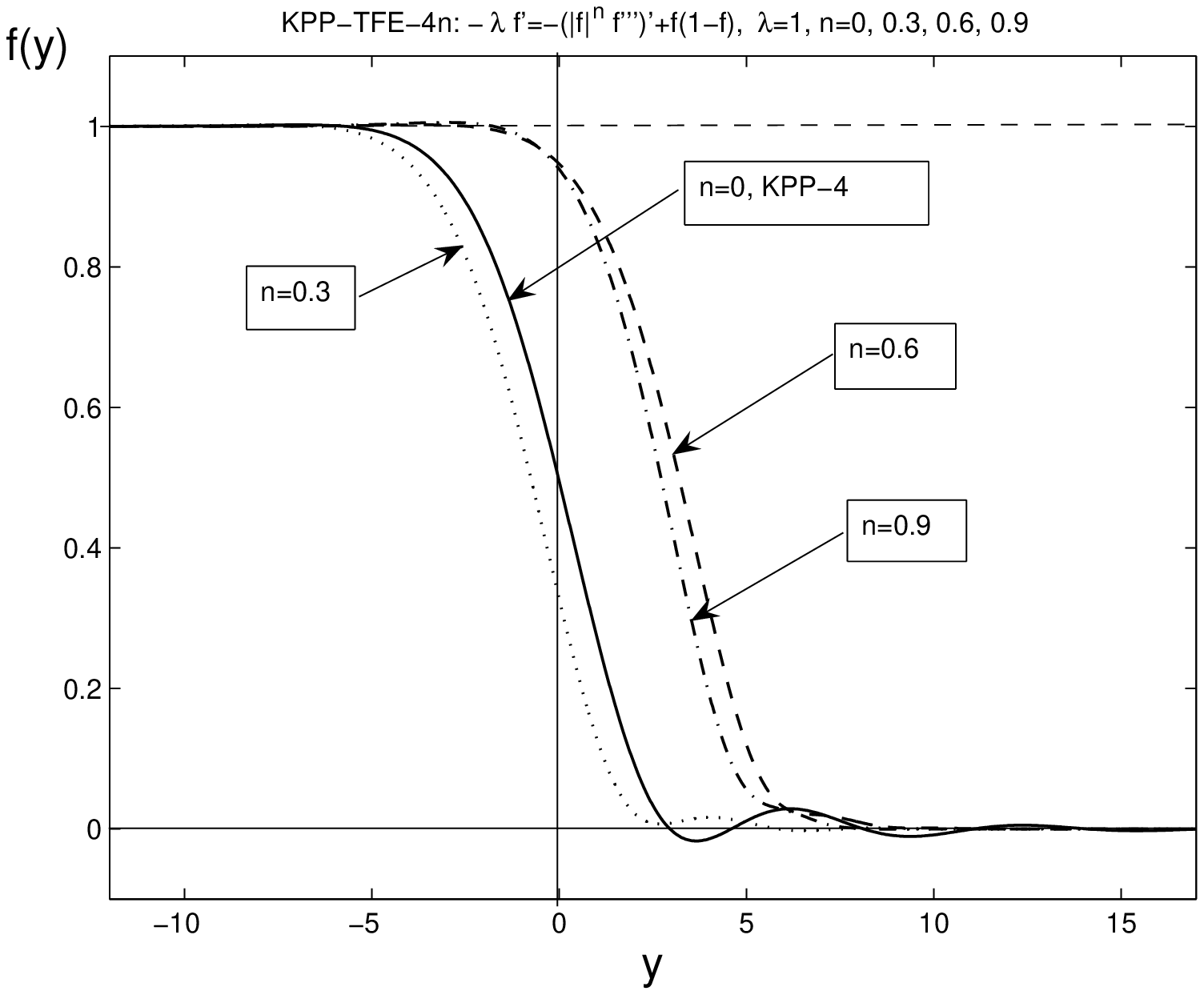}  
\vskip -.3cm
  \caption{The TW profiles $f(y)$ of
(\ref{tfe1}), \ef{BC1}  for
 $\l=1$, $n=0, \,0.3,\, 0.6$, and $0.9$.}
 \label{Ftfe1}
\end{figure}

 We next briefly review some results for the KPP problem for the parabolic {\em thin film equation}
 (TFE--4)
  ($n>0$)
  \be
  \label{tfe1}
  u_t=-(|u|^n u_{xxx})_x+u(1-u) \LongA -\l f'=-(|f|^n
  f''')'+f(1-f).
   \ee
   In Figure \ref{Ftfe1}, we show TW profiles for $\l=1$ and
   sufficiently small $n=0.3$, 0.6, and 0.9, and, for comparison,
   we present the ``linear" one for $n=0$ (i.e., for the KPP--4
   problem in \cite{GKPPI}). Next Figure \ref{Ftfe2} shows
   profiles for larger $n=1.5$, 2, 3, and 4.

In Figure \ref{Ftfe3}, we show an enlarged image of the behaviour
of $f(y)$ from Figure \ref{Ftfe2} for large $y \in [10,17]$. It is
seen that $f(y)$ becomes less oscillatory close to interfaces as
$n$ reaches about 2. We believe that a changing sign feature for
$n \ge 2$ in this figure can be related to the necessary
regularization of the degenerate thin film operator in \ef{tfe1},
where we replace
 \be
 \label{tfe4}
 |f|^n \mapsto (\e^2+f^2)^{\frac n2}, \quad \mbox{with} \quad
 \e=10^{-3},
  \ee
  since smaller $\e=10^{-4}$, etc., always led to divergence of
  the numerical scheme; see below.

    Note that the
   convergence for larger $n$ got very slow; for $n=4.3$, we
   still got some profile, which is close to that for $n=4$ in Figure \ref{Ftfe2},
    while, for $n \ge 4.4$, ``a singular
   Jacobian" (a full non-convergence) always appeared.

The oscillatory behaviour of $f(y)$ near the interface was studied
in \cite{Gl4, GHarCentre}. It was shown that a periodic
oscillatory component $\var(s)$ (a full analogy of that in
\ef{le3}, \ef{le4}) exists up to a critical {\em homoclinic
bifurcation} exponent $n_{\rm h}$,
\be
  \label{pp1}
   \mbox{$
  0<n<n_{\rm h} \in (\frac 32,n_{\rm +}), \quad \mbox{where}
   \quad n_{\rm +}= \frac 9{3+\sqrt 3}=1.9019238... \, .
    $}
    \ee
Numerically, $n_{\rm h}$ is given by
  \be
\label{n**1}  n_{\rm h}= 1.75987... \, .
  \ee
 Therefore, it is expected that, for $n>n_{\rm h}$, $f(y)$
 exhibits either a finite number of zeros close to the interface,
 or even becomes nonnegative (at least, for $n\ge 2$).


 \begin{figure}
\centering
\includegraphics[scale=0.85]{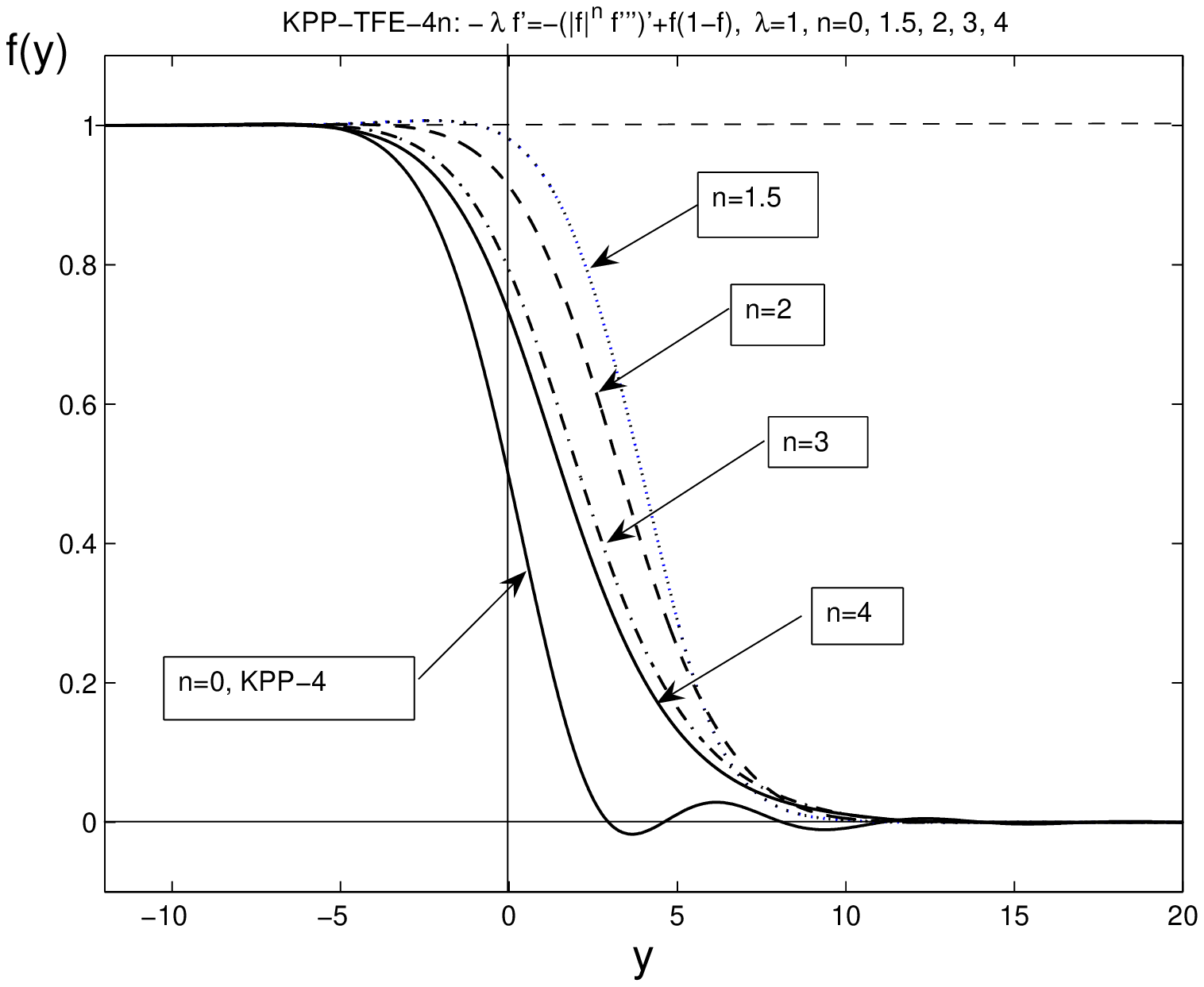}  
\vskip -.3cm
  \caption{The TW profiles $f(y)$ of
(\ref{tfe1}), \ef{BC1}  for
 $\l=1$, $n=1.5, \,2,\,3$, and $4$.}
 \label{Ftfe2}
\end{figure}



 \begin{figure}
\centering
\includegraphics[scale=0.85]{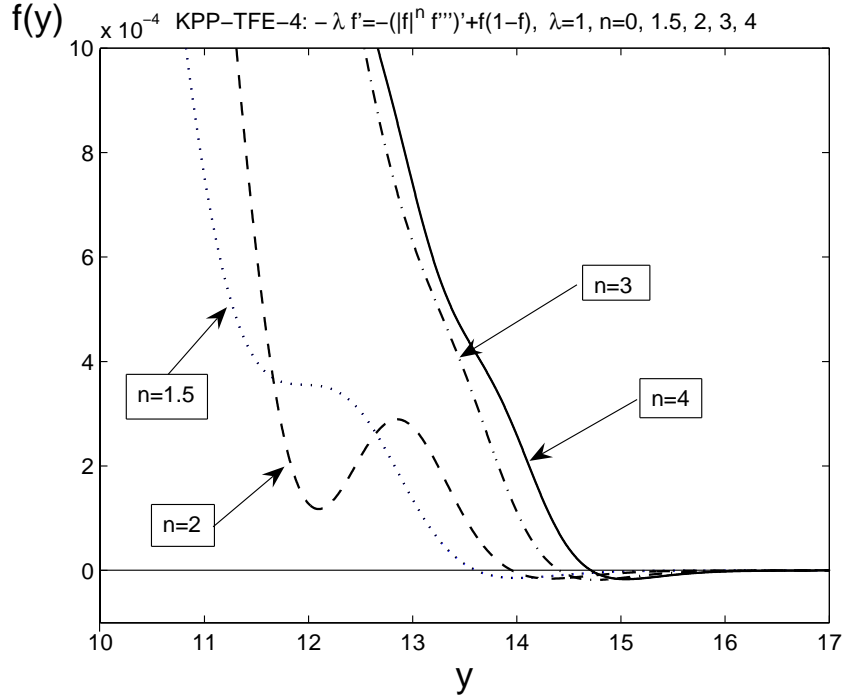}  
\vskip -.3cm
  \caption{Enlarged behaviour of  TW profiles $f(y)$ close to zero of
(\ref{tfe1}), \ef{BC1} from Figure \ref{Ftfe2};
 $\l=1$, $n=1.5, \,2,\,3$, and $4$.}
 \label{Ftfe3}
\end{figure}


Finally, existence of a $\log t$-shift of moving front for the PDE
KPP--TFE--4 problem \ef{tfe1} is shown precisely in the same way
as in Section \ref{SDiscr}.

\section{The quasilinear KPP--(6,1)$n$ problem}
\label{S3}

We now, more briefly than above, consider the KPP--(6,1)$n$
problem \ef{PME6} and its ODE counterpart
 \be
 \label{E7}
  - \l f'=(|f|^nf)^{(6)}+f(1-f), \quad F=|f|^n f,
  \ee
with singular  boundary conditions \ef{BC1}.

\begin{figure}
\centering
\includegraphics[scale=0.85]{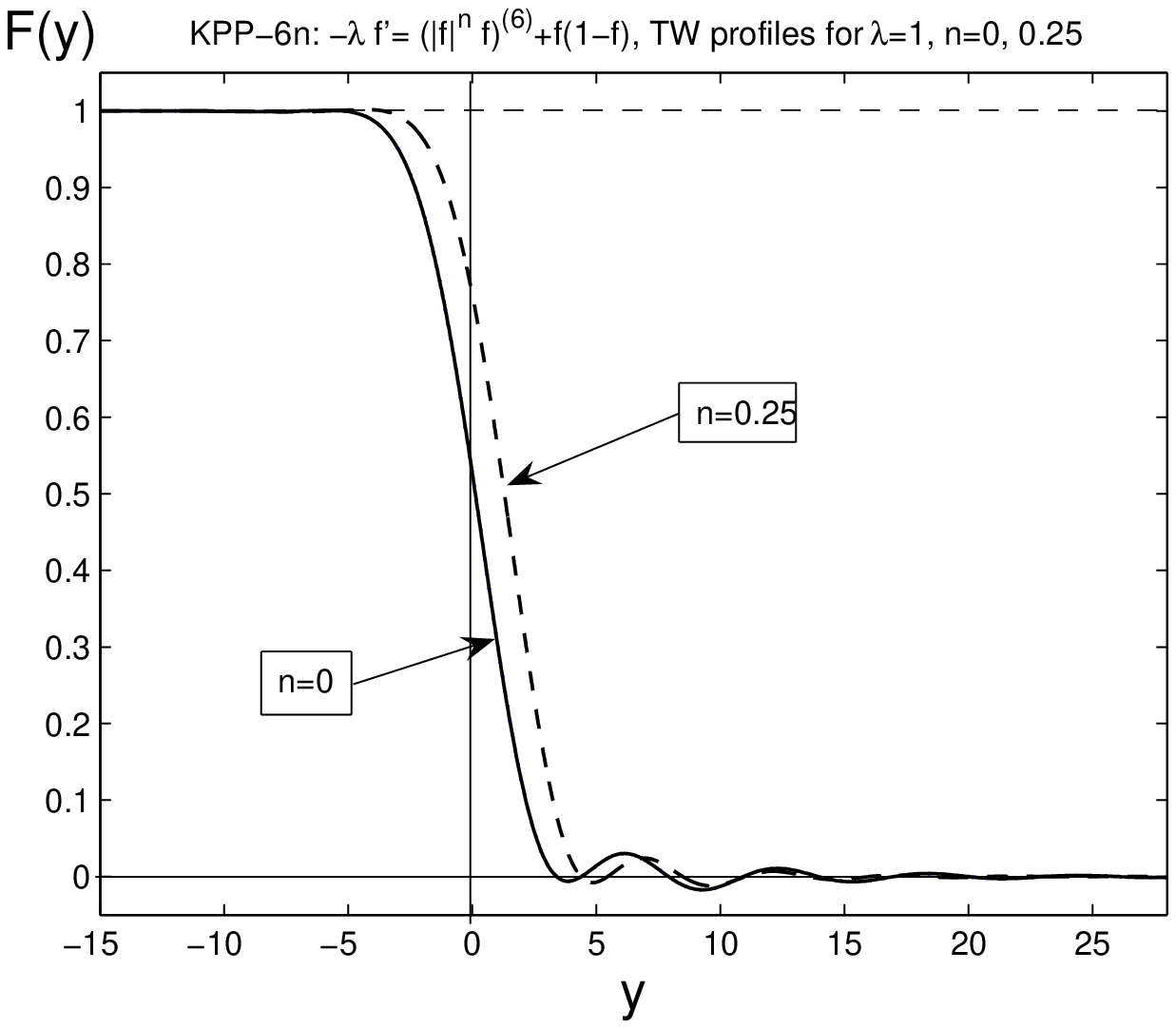}  
\vskip -.3cm
  \caption{The TW profiles $F(y)$ satisfying
(\ref{FF11}), \ef{BC1} for
 $\l=1$,  and $n=0.25$, $n=0$.}
 \label{Fn61}
\end{figure}

Figure \ref{Fn61} shows TW profiles for $\l=1$ and $n=0.25$. For
comparison, we also put the ``semilinear" profile for $n=0$. This
confirms a rather non-surprising fact that there exists a
continuous dependence of $f(y)$ on $n \to 0^+$ (this even can be
proved for such ODEs).

\begin{figure}
\centering
\includegraphics[scale=0.85]{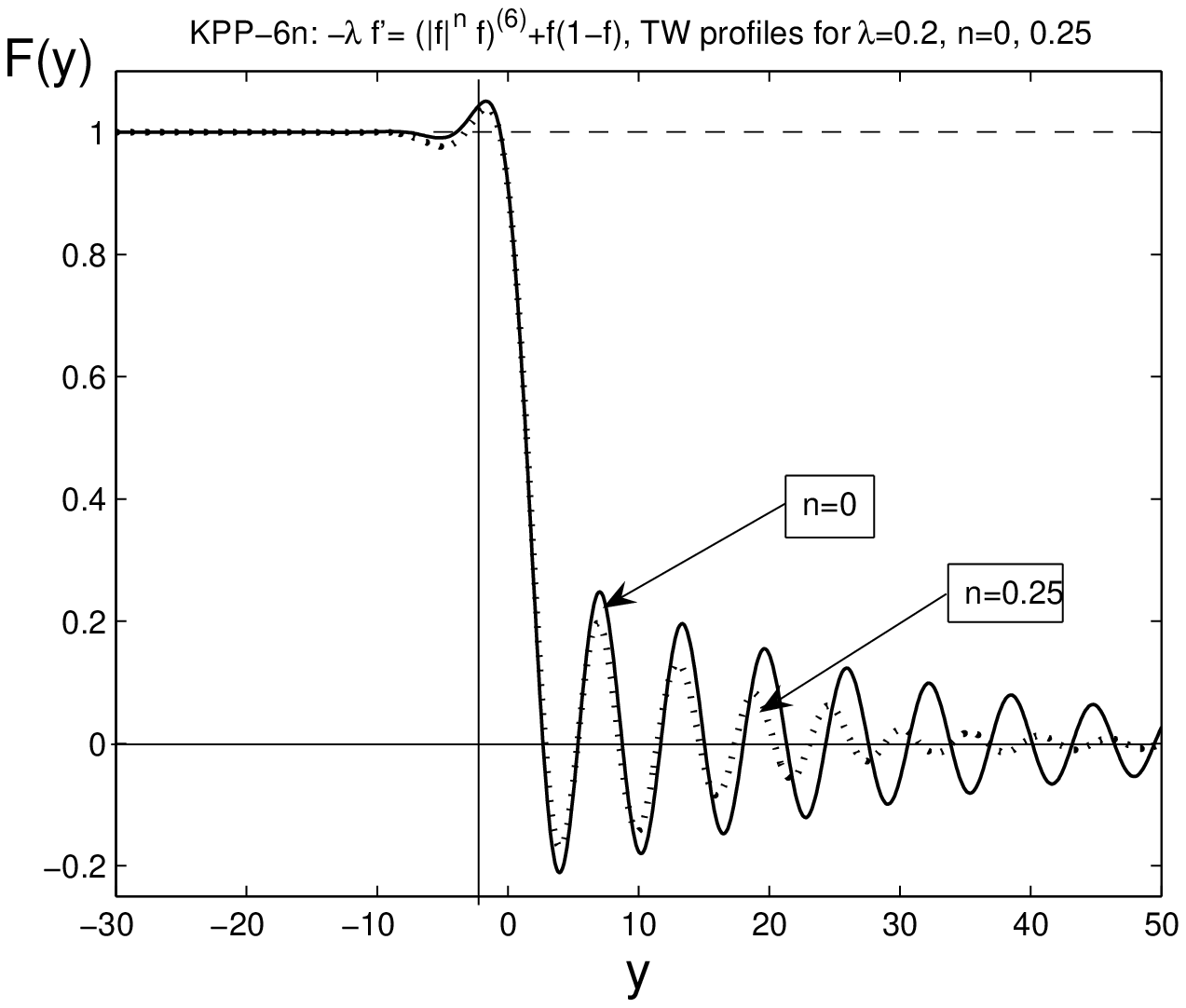}  
\vskip -.3cm
  \caption{The TW profiles $F(y)$ satisfying
(\ref{FF11}), \ef{BC1} for
 $\l=0.2$,  and $n=0.25$, $n=0$.}
 \label{Fn62}
\end{figure}

In Figure \ref{Fn62}, we present a similar comparison for the more
oscillatory case  $\l=0.2$, where we again observe existence of  a
clear ``homotopy" deformation as $n \to 0$.


\ssk

A proper dimensional analysis of the linearized bundle as $y \to
-\iy$ (i.e., as $f \to 1$) is performed similarly, as in
\cite{GKPPI}.

The oscillatory periodic-like behaviour at the finite interface as
$y \to y_0^-$ (i.e., as $f \to 0$) is performed as in \cite{GBl6};
see also \cite[p.~142]{GSVR}.

Note  that the corresponding blow-up problem for $n \in (0,1)$
(cf. the ODE \ef{t13})
 \be
 \label{bl1}
 (|f|^n f)^{(6)}=f^2
 \ee
 is now more difficult since \ef{bl1} admits oscillatory
 solutions, which can be studied as in \cite[\S~7]{Gl4} and in \cite{GBl6} by
 introducing an {\em oscillatory component} represented by
 periodic or other functions.

 \section{Very briefly on KPP--8$n$}
\label{S810}

\begin{figure}
\centering
\includegraphics[scale=0.85]{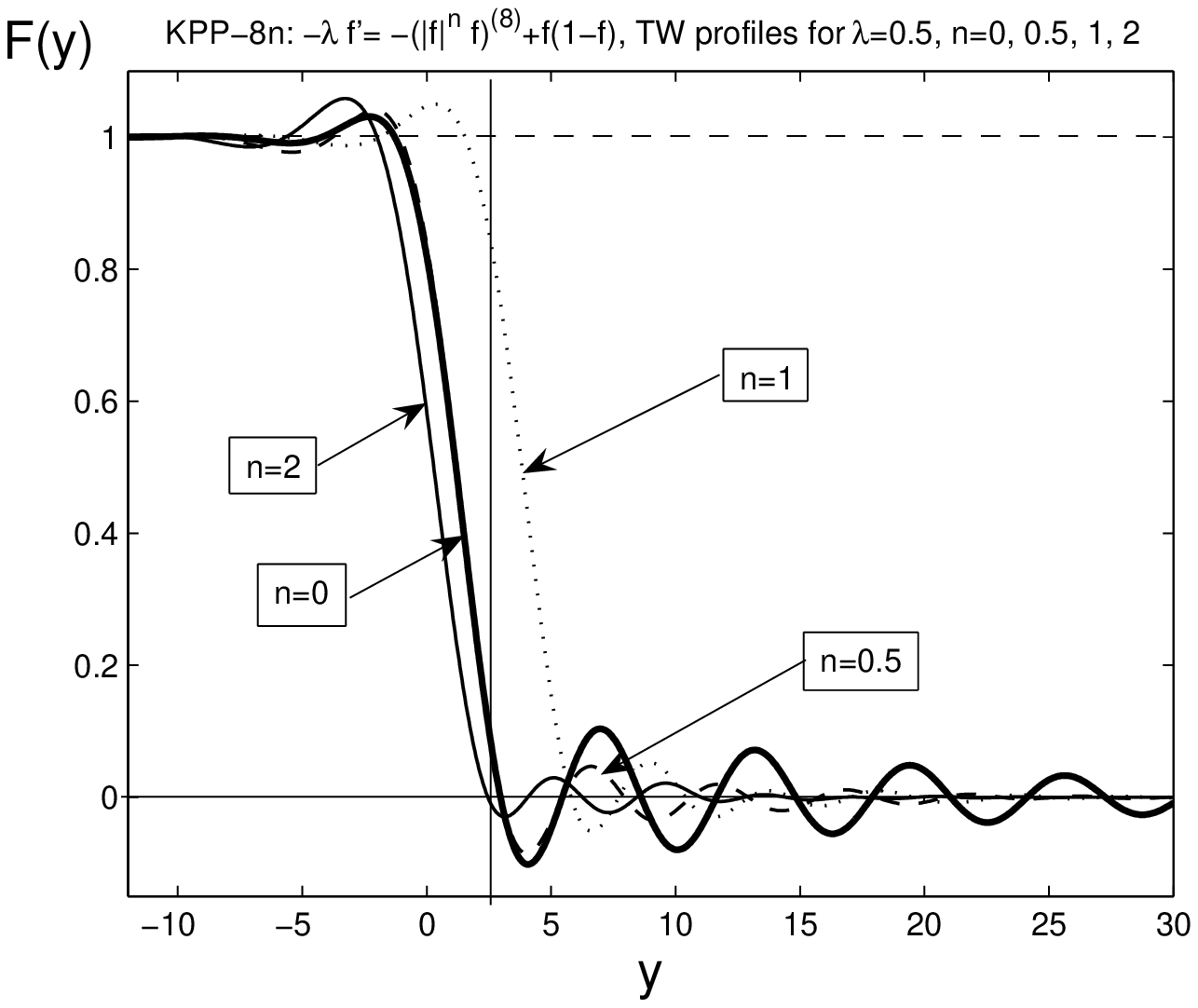}  
\vskip -.3cm
  \caption{The TW profiles $f(y)$ satisfying
(\ref{8.1}), \ef{BC1} for
 $\l=0.5$ and $n=0, \, 0.5, 1, \, 2$.}
 \label{F81n}
\end{figure}

 In the KPP--8$n$, we deal with the following ODE problem:
  \be
  \label{8.1}
  u_t=-D_x^8(|u|^n u)+u(1-u) \LongA (\mbox{ODE}) \quad -
 \l f'=-(|f|^n f)^{(8)}+f(1-f).
  \ee
  Figure \ref{F81n} shows  TW profiles $F$ for $\l=0.5$ and $n=0.5, \, 1, \,2$, plus,
   for comparison, for
  $n=0$ from \cite[\S~4]{GKPPI}.
 This shows a principal positive answer on the TW existence question.
It is seen that oscillations about $f=0$ decrease as $n$
increases, though finite interfaces for $n>0$ (unlike $n=0$ with
an exponential decay as $y \to +\iy$) are obviously invisible.


Stable bundles as $f \to 1$ ($y \to -\iy$) are studied as in
\cite{GKPPI}, while, close to finite interfaces, as $f \to 0$,
periodic oscillatory components  for seventh-order ODEs are
constructed as in \cite[p.~142]{GSVR}.

\section{The origin of $\log t$-shift: centre subspace balancing}

 \label{SDiscr}

 We now explain the origin of possible $\log t$-shiftings (say, retarding) from the
TW for the PDE higher-order quasilinear parabolic  KPP--problems
\ef{PME4} or \ef{PME6}.
Note that the actual proof of such $\log t$-drift assumes
 a delicate matching of the solution behaviour on compact subsets in the TW
$y$-variable, i.e., in the {\em Inner  Region}, with a remote {\em
Outer Region} for $y \gg 1$, where the influence of the nonlinear
term $-u^2$ is negligible, and the actual behaviour is governed by
the quasilinear bi- (or tri-) harmonic operator. Such a matching,
eventually, is supposed to describe classes of initial data, for
which such ($\pm$) $\log t$-shifting (or no shifting at all, say,
a pure convergence to a TW in its moving frame).
We do not present here this kind of a procedure of {\em matching}
of asymptotic expansions of those two regions, which is extremely
difficult in the  quasilinear/degenerate case. Therefore, we
restrict to an Inner expansion analysis.

\subsection{Linearization and rescaled equation}

Thus, we consider a KPP-type problem for a quasilinear PDE
 \be
 \label{K1}
  u_t=\AAA(u)+u(1-u) \inB \re \times \re_+,
  \ee
  with some proper step-like data  $u_0(x)$, with, however, some
  positive, negative, or even oscillatory ``tails" as $y \to
  +\infty$, in order to create a necessary $\pm \log t$-shift.
  Here, in \ef{K1}, $\AAA(u)$ is a proper homogeneous isotropic quasilinear differential
  operator satisfying some extra conditions specified below.
 We  fix, as key examples, the quasilinear
  bi- or tri-Laplacian operators of the {\em porous medium type}
   \be
   \label{Lapl22}
    \AAA(u)=-D_x^4(|u|^n u) \orA   \AAA(u)=D_x^6(|u|^n u), \quad n>0.
   \ee

   We assume that, for some fixed $\l_0>0$, the corresponding ODE problem
   \be
   \label{K2}
   -\l_0 f'= \AAA(f)+f(1-f),
   \ee
   with the conditions \ef{BC1} admits a unique solution $f$.

   Attaching the solution $u(x,t)$ to the
    front moving and setting, for convenience, $x_f(t) \equiv \l_0t -g(t)$,
    the PDE reads
    \be
    \label{K21}
    u(x,t)=v(y,t), \quad y=x- \l_0 t+g(t) \LongA v_t= \AAA(v) + v(1-v)
    + \l_0 v_y - g'(t) v_y.
     \ee
   We next linearize \ef{K21} by setting
    \be
    \label{K3}
    v(y,t)=f(y)+w(y,t),
     \ee
     that yields the following perturbed equation:
      \be
      \label{gg2}
       \begin{aligned}
      w_t=  &\BB w -g'(t) f' - g'(t) w_y -w^2,\quad \mbox{where} \ssk\ssk \\
    &  \BB w= \AAA'(f) w + (1-2f)w
       + \l_0 w_y \andA \AAA'(f)w=(n+1)\AAA(|f|^n w).
        \end{aligned}
       \ee
 Assuming that, in this $g(t)$-moving frame, there exists the
 convergence as in \ef{TW3}, so that $w(t) \to 0$ as $t \to +\iy$,
 one can see that the leading non-autonomous perturbation
in \ef{gg2} is the second term on the right-hand side. However, as
we show, the last two terms, though negligible, will define a
proper $\log t$-shift of the front.

\subsection{Where $\log t$-shift comes from}

Again, as in the semilinear case $n=0$, we note that the rescaled
equation \ef{gg2} is essentially {\em non-autonomous} in time, so
we cannot use powerful  tools of nonlinear semigroup theory; see
\cite{Lun}. However, using  a formal asymptotic approach, we will
trace out some definite {\em centre subspace} behaviour after an
extra rescaling and balancing of non-autonomous perturbations.

 Thus, as usual (see Introduction), we assume that $g'(t) \to 0$ as
 $t \to +\iy$ sufficiently fast, i.e., at least algebraically, so
 that
  \be
  \label{gg1}
  |g''(t)| \ll |g'(t)| \forA t \gg 1.
  \ee
Under the hypothesis \ef{gg1}, the only possible way to balance
{\em all} the terms therein (including the quadratic one $-w^2$)
for $t \gg 1$ is to assume the asymptotic separation of variables:
 \be
 \label{gg3}
 w(y,t) = g'(t) \psi(y)+ \e(t) \var(y)+... \whereA
 |\e(t)|
 \ll |g'(t)| \asA t \to + \iy.
  \ee
 Here, we omit higher-order perturbations. Substituting \ef{gg3}
  into \ef{gg3} yields
   \be
   \label{gg31}
    \begin{aligned}
    & g''(t) \psi + \e'(t) \var+...
   =  \,g'(t)(\BB \psi-f') \ssk \\
   & + \e(t)\BB
   \var -(g'(t))^2(\psi'+\psi^2)-g'(t)\e(t)(\var'+2\var \psi)-\e^2(t)
   \var^2+...\, .
 \end{aligned}
   \ee

Using \ef{gg1} and \ef{gg3} in
 balancing first the leading terms of the order
  $O(g'(t))$ yields the elliptic equation for $\psi$:
   \be
   \label{gg4}
    \tex{
 O(g'(t))\,\,\big(=O(\frac 1t)\big): \quad  \BB \psi - f'=0.
   }
   \ee
   Then balancing the rest of the terms in \ef{gg31} requires
   their asymptotic equivalence,
 \be
 \label{hh41}
  \tex{
 g''(t) \sim -(g'(t))^2 \sim \e(t), \,\,\, \mbox{i.e.,} \,\,\,
  g(t)= k \log t, \,\, \mbox{$g'(t)= \frac kt, \,\,
 g''(t)=- \frac k{t^2}$,} \,\, \e(t)= \frac 1{t^2}.
 }
 \ee
 Then, we obtain the second inhomogeneous singular
 Sturm--Liouville problem for $\var$:
  \be
  \label{SL1}
   \tex{
  O\big(\frac 1{t^2}\big): \quad
  \BB \var= k \psi + k^2(\psi'+\psi^2).
 }
   \ee

    Thus, the first simple asymptotic ODE in \ef{hh41} gives the $\log
    t$-dependence as in \ef{1.3}. Finally, we arrive at the
    following system for $\{\psi,\var\}$:
     \be
     \label{sys21}
     \left\{
     \begin{aligned}
& \BB \psi = f', \\
   &
   \BB \var= k \psi + k^2(\psi'+\psi^2).
    \end{aligned}
     \right.
 \ee
  Solving this system, with typical boundary conditions as in
  \ef{BC1}, allows then continue the expansion of the solutions of
  \ef{gg2}
  close to an  ``affine ({\rm i.e., shifted via $f'$ on the RHS}) centre subspace" of $\BB$
  governed by the spectral pair obtained by translation in
  \ef{K2}:
   \be
   \label{pair1}
   \hat \l_0=0 \andA \hat \psi_0(y)= f'(y).
    \ee
 The asymptotic expansion for $t \gg 1 $ then takes the form
  \be
  \label{as22}
  \tex{
  w(y,t)= \frac kt \, \psi(y) + \frac 1{t^2} \,\var(y) + ...\, ,
  }
  \ee
  which can be easily extended by introducing further terms, with
  similar inhomogeneous Sturm--Liouville problems for the
  expansion coefficients.

As for $n=0$, $\BB$ does not have a discrete spectrum, so we
cannot get a simple algebraic equation for $k$ by demanding the
standard orthogonality of the right-hand side in the second
equation in \ef{sys21} to the adjoint eigenvector $\hat \psi^*_0$
of $\BB^*$ in some
``weighted" metric $\langle \cdot,\cdot
\rangle_*$ (in which the adjoint operator $\BB^*$ is obtained, if
any), like
 \be
 \label{as23}
k: \quad   \langle k \psi + k^2(\psi'+\psi^2), \, \hat \psi_0^*
\rangle_* =0.
 \ee
 Therefore, it seems, the system \ef{sys21} cannot itself  determine
    the actual value of $k$ therein. As we have mentioned, the latter requires a
    difficult matching analysis in Inner and Outer Regions, which,
    for the KPP--4$n$ (and all other problems) remains an open
    problem.






\begin{thebibliography}{10}






 \bibitem
  {And10}
 F.~El Adnani and H.~Talibi Alaoui, {\em Travelling front
 solutions in nonlinear duffusion degenerate Fischer-KPP anfd
 Nagumo equations via the Conley index}, Topol. Math. Nonl.
 Anal., {\bf 35} (2010), 43--60.


 \bibitem
  {Ar80}
 D.G.~Aronson, {\em Density-dependent inteaction-diffusion systems
 in dynamics and modelling of reactive systems}, In: Dynamics and Modelling of Reactive Systems,
 W.E.~Stewart {et al} Eds, Acad. Press, New York,
 1980, 161--176.

 \bibitem
  {Atk81}
  C.~Atkinson, G.E.H.~Reuter, and C.J.~Ridler-Rowe, {\em
  Travelling wave solutions  for some nonlinear diffusion
  equations}, SIAM. J.~Math. Anal., {\bf 12} (1981), 880--892.



\bibitem{B}
  G.I. Barenblatt, {Similarity, Self-Similarity,
Intermediate Asymptotics}, Consultant Bureau, New York, 1978.

\bibitem
 {Bern88}
 F.~Bernis, \emph{Source-type solutions of fourth order
 degenerate parabolic equations}, In: {Proc. Microprogram
 Nonlinear Diffusion Equation and Their Equilibrium States},
  W.-M.~Ni, L.A.~Peletier, and J.~Serrin, Eds., MSRI Publ., Berkeley,
 California, Vol. 1, New York, 1988, pp. 123--146.



\bibitem 
 {BF1}
 F.~Bernis and A.~Friedman, \emph{Higher order nonlinear degenerate
 parabolic equations}, J.~Differ. Equat., \textbf{83} (1990), 179--206.







\bibitem 
 {BMc91}
 F.~Bernis and J.B.~McLeod, \emph{Similarity solutions of a higher order
nonlinear
 diffusion
 equation}, Nonl. Anal., TMA, \textbf{17} (1991), 1039--1068.







\bibitem{Bert01}
 A.L.~Bertozzi, A.~M\"unch, M.~Shearer, and K.~Zumbrun, {\em
Stability of compressive and undercompressive thin film travelling
waves},
 Euro J.~Appl. Math., {\bf 12} (2001),
253--291.

\bibitem{Bert99}
 A.L.~Bertozzi and M.~Shearer, {\em Existence of undercompressing
 travelling waves in thin film equations},
 SIAM J.~Math. Anal., {\bf 32} (1999),
194--213.


\bibitem{Br}
 M.~Bramson, {\em Convergence of solutions of the
Kolmogorov equation to travelling waves}, Memoirs of  Amer. Math.
Soc., {\bf 44} (1983), 1--190.



  \bibitem{Collet90}
   P.~Collet and J.-P.~Eckmann, {Instabilities and Fronts in
   Extended Systems}, Pinceton Univ. Press, Princeton, NJ, 1990.






\bibitem
{Deim} K. Deimling, {\rm Nonlinear Functional Analysis},
Springer-Verlag, Berlin/Tokyo, 1985.








\bibitem 
{EidSys}
 S.D. Eidelman, {Parabolic Systems,}  North-Holland Publ.
Comp., Amsterdam/London, 1969.




   \bibitem{Gl4}
J.D.~Evans, V.A.~Galaktionov, and J.R.~King, 
\emph{Source-type solutions of the fourth-order unstable thin film
equation}, Euro J.~Appl. Math., {\bf 18} (2007), 273--321.



      \bibitem
      {GBl6}
J.D.~Evans, V.A.~Galaktionov, and J.R.~King, {\em Unstable
sixth-order thin film equation. I. Blow-up similarity solutions;
II. Global similarity patterns},
 {Nonlinearity}, {\bf 20} (2007), 1799--1841, 1843--1881.


 \bibitem
  {EGW1}
J.D.~Evans, V.A.~Galaktionov, and J.F.~Williams, {\em Blow-up and
global asymptotics of  the limit unstable Cahn-Hilliard equation},
SIAM J. Math. Anal., {\bf 38} (2006), 64--102.


\bibitem
 {Fife77}
  P.~Fife and J.B.~McLeod, {\em The approach of solutions
 of nonlinear diffusion equations to travelling front solutions},
 Arch. Rat. Mech. Anal., {\bf 65} (1977), 335--361.


\bibitem  
{GalTW81}  V.A.~Galaktionov, {\em  Some properties of travelling
waves in a medium with non-linear thermal conductivity and a
source of heat,} {USSR Comput. Math. Math. Phys.,} {\bf 21}
(1981), 167--176.




\bibitem 
{Gal2m}
   V.A.~Galaktionov, {\em On a spectrum of blow-up patterns
for a higher-order semilinear parabolic equation}, Proc. Royal
Soc. London A, {\bf 457} (2001), 1--21.





\bibitem
 {GalCr}
  V.A. Galaktionov, {\em Critical global asymptotics in
  higher-order semilinear parabolic equations}, Int.
  J. Math. Math. Sci., {\bf 60} (2003),  3809--3825.


\bibitem 
 {GalGeom}
 V.A.~Galaktionov, {\rm Geometric Sturmian  Theory of Nonlinear
 Parabolic Equations and Applications}, Chapman and Hall/CRC, Boca Raton,
Florida,
 2004.



  \bibitem
 {GalRDE4n}
 V.A.~Galaktionov,
 {\em Countable branching of similarity solutions of higher-order
 porous medium type equations}, Adv. Differ. Equat., {\bf 13}
 (2008), 641--680.



\bibitem
 {GKPPI}
 V.A.~Galaktionov, {\em Towards the KPP--problem  and
  ${{\log t}}$--front  shift for higher-order nonlinear PDEs
I. Bi-harmonic and  other parabolic equations}, to appear
(arXiv:1210.3513).



\bibitem
 {GKPPIII}
 V.A.~Galaktionov, {\em Towards the KPP--problem and
  ${{\log t}}$--front  shift for higher-order nonlinear PDEs III.
 Dispersion and hyperbolic equations},
 to appear (available in arXiv.org).

\bibitem
{GHarCentre}
 V.A.~Galaktionov and P.J.~Harwin,
 {\em On centre subspace behaviour in thin film
equations}, SIAM J.~Appl. Math., {\bf 69} (2009), 1334--1358.














\bibitem
 {GMPSobI}
 V.A.~Galaktionov, E.~Mitidieri,  and S.I.~Pohozaev,
 \emph{Variational approach to complicated similarity solutions
of higher-order
 nonlinear evolution equations of parabolic, hyperbolic, and nonlinear dispersion
 types},
In: Sobolev Spaces in Mathematics. II, Appl. Anal. and Part.
Differ. Equat., Series: Int. Math. Ser., Vol. {\bf 9}, V.~Maz'ya
Ed., Springer, 2009 (an earlier preprint: arXiv:0902.1425).

\bibitem
 {GMPSobII}
 V.A.~Galaktionov, E.~Mitidieri,  and S.I.~Pohozaev,
 {\em Variational approach to complicated similarity solutions
of higher-order
 nonlinear  PDEs. II},  Nonl. Anal.: RWA, {\bf 12} (2011),
 2435--2466
  (arXiv:1103.2643).




\bibitem
 {GSVR} V.A.~Galaktionov and S.R.~Svirshchevskii, Exact Solutions and
 Invariant Subspaces of Nonlinear Partial Differential Equations in Mechanics and Physics,
  Chapman$\,\&\,$Hall/CRC, Boca Raton,
Florida,
 2007.





\bibitem{Ga}
 J. G\"artner, {\em Location of wave front for the
multidimensional K-P-P equation and  brownian first exit
densities},  Math. Nachr., {\bf 105} (1982),  317--351.


\bibitem
 {Gild04}
 B.H.~Gilding and R.~Kersner, {\rm
 Travelling Waves in Nonlinear
 Diffusion-Convection Reaction}, Birkh\"auser Verlag, Basel, 2004.


 \bibitem
 {Gild05}
 B.H.~Gilding and R.~Kersner, {\em A Fisher/KPP-type equation with
 density-depemdent diffusion and convection: travelling-wave
 solutions}, J.~Phys. A, {\bf 38} (2005), 3367--3379.

\bibitem
 {Grin87}
   P.~Grinrod and B.G.~Sleeman, {\em Weak
 travelling fronts for population models with density-dependent
 dispersion},
Math. Meth. Appl. Sci.,
   {\bf 9} (1987), 576--586.

\bibitem
 {Had75}
   K.P.~Hadeler and F.~Rothe, {\em
 Travelling fronts in non-linear diffusion equations}, J.~Math.
 Biol.,
   {\bf 2} (1975), 251--263.

\bibitem
 {Hos86}
   Y.~Hosono, {\em
 Travelling wave solutions for some density dependent diffusion
 equations}, Japan. J.~Appl. Math.,
   {\bf 3} (1986), 163--196.




\bibitem
 {Hu2008}
   C.~Hu, {\em Stability of under-compressive waves with second
   and fourth order diffusion}, Discr. Cont. Dyn. Syst.,
   {\bf 22} (2008), 629--662.





\bibitem
 {KKVV00}
  W.D.~Kalies, J.~Kwapisz, J.B.~VandenBerg, and
  R.C.A.M.~VanderVorst, {\em Homotopy classes for stable periodic
  and chaotic patterns in fourth-order Hamiltonian systems,}
  Commun. Math. Phys., {\bf 214} (2000), 573--592.

\bibitem
 {King03}
J.R.~King and P.M.~McCare, {\em On the Fischer--KPP equation with
fast nonlinear diffusion}, Proc. Roy. Soc., Sect.~A, {\bf 459}
(2003), 2529--2546.






\bibitem  
{KPP}   A.N.~Kolmogorov, I.G.~Petrovskii, and  N.S.~Piskunov, {\em
 Study of the diffusion equation with  growth of the
quantity of matter and its application to a biological problem,}
{Byull. Moskov. Gos. Univ., Sect. A,} {\bf 1} (1937), 1--26.
English. transl. In: {Dynamics of Curved Fronts,} P.~Pelc\'e, Ed.,
Acad. Press, Inc., New York, 1988, pp.~105--130.


\bibitem{KrasZ}
M.A.~Krasnosel'skii and P.P.~Zabreiko, {Geometrical Methods of
Nonlinear
  Analysis}, Springer-Verlag, Berlin/Tokyo, 1984.






\bibitem
{Lun} A. Lunardi, { Analytic Semigroups and Optimal Regularity in
Parabolic Problems}, Birkh\"auser, Basel/Berlin, 1995.

\bibitem
 {Man2010}
 M.B.A.~Mansour, {\em Travelling wave solutions for the extended
 Fisher/KPP equation},
  Reports Math. Phys., {\bf 66} (2010), 375--383.

 \bibitem
 {PaVa91}
   A.~De Pablo and J.L.~Vazquez, {\em
 Travelling waves and finite propagation in a reaction-diffusion
 equation},
J.~Differ. Equat.,
   {\bf 93} (1991), 19--61.





\bibitem
{PelTroy} L.A.~Peletier and W.C.~Troy, {\rm Spatial Patterns.
Higher Order Models in Physics and Mechanics}, Birkh\"auser,
Boston/Berlin, 2001.





 \bibitem  
{SGKM}    A.A.~Samarskii, V.A.~Galaktionov, S.P.~Kurdyumov, and
A.P.~Mikhailov, {Blow-up in Quasilinear Parabolic Equations,}
 Walter de
Gruyter \& Co., Berlin, 1995.

\bibitem
 {San95}
   F.~Sanchez-Garduno and P.K.~Maini, {\em
 Travelling wave phenomena in some degenerate reaction-diffusion equations},
  J.~Differ.
 Equat.,
   {\bf 117} (1995), 1--41.

   \bibitem
 {San97}
   F.~Sanchez-Garduno and P.K.~Maini, {\em
 Travelling wave phenomena in non-linear diffusion degenerate Nagumo equations},
  J.~Math.
 Biol.,
   {\bf 35} (1997), 731--728.








 \bibitem
 {VV02}
   J.B.~Van Den Berg and
  R.C.~Vandervorst, {\em Stable
   patterns for fourth-order parabolic equations,}
  Duke Math.~J., {\bf 115} (2002), 513--558.


\bibitem
 {Vaz07}
  J.L.~Vazquez, {\em Porous medium flow in a tube. Travelling waves
  and KPP behaviour}, Commun. Contemp. Math., {\bf 9} (2007),
  731--751.







\end{thebibliography}
\end{document}